\documentclass[final,leqno]{siamltex}

\usepackage{amsfonts}
\usepackage{amsmath}
\usepackage{graphicx}
\usepackage{graphics}
\usepackage{subfigure}
\usepackage{caption}

\title{Entropy Conservative Schemes and the receding flow problem\thanks{This work was
supported by the AFOSR through grant no. FA9550-16-1-030.}} 

\author{
  Ayoub GOUASMI\thanks{University of Michigan, Dept. of Aerospace Engineering, Ann Arbor, MI
    (gouasmia@umich.edu).}
  \and
    Scott Murman\thanks{NASA Ames Research Center, NASA Advanced Supercomputing Division, Moutain View, CA (scott.m.murman@nasa.gov).}
  \and
  Karthik Duraisamy\thanks{University of Michigan, Dept. of Aerospace Engineering, Ann Arbor, MI
    (kdur@umich.edu).}
}

\begin{document}
\maketitle

\begin{abstract}
This work delves into the family of entropy conservative (EC) schemes introduced by Tadmor. The discussion is centered around the Euler equations of fluid mechanics and the receding flow problem extensively studied by Liou. This work is motivated by Liou's recent findings that an abnormal spike in temperature observed with finite-volume schemes is linked to a spurious entropy rise, and that it can be prevented in principle by conserving entropy. While a semi-discrete analysis suggests EC schemes are a good fit, a fully discrete analysis based on Tadmor's framework shows the non-negligible impact of time-integration on the solution behavior.  An EC time-integration scheme is developed to show that enforcing conservation of entropy at the fully discrete level does not necessarily guarantee a well-behaved solution. 
\end{abstract}

%\begin{keywords}
%entropy conservation, time-integration, fully discrete schemes, conservation %laws
%\end{keywords}

%\begin{AMS}
%35B06, 35L65, 35Q31, 35Q35
%\end{AMS}

%\pagestyle{myheadings}
\thispagestyle{plain}
\markboth{AYOUB GOUASMI, SCOTT MURMAN, AND KARTHIK DURAISAMY}{ENTROPY CONSERVATIVE SCHEMES AND THE RECEDING FLOW PROBLEM}

\section{Introduction}

In this paper, we consider general hyperbolic systems of conservation laws of the form
\begin{equation}\label{eq:base}
\frac{\partial u}{\partial t} + \frac{\partial f}{\partial x} = 0, \ u = u(x,t) \in \mathbb{R}^{\mathbb{N}}, \ x \in \mathbb{R}, t > 0,
\end{equation}
and their numerical solution using finite-volume schemes. In equation (\ref{eq:base}), $f$ is a smooth nonlinear flux function of the unknown $u$. We also assume that the system (\ref{eq:base}) has a convex extension in the sense that an additional conservation law \cite{Lax}: 
\begin{equation}\label{eq:base_entropy}
\frac{\partial U}{\partial t} + \frac{\partial F}{\partial x} = 0,
\end{equation}
where $U = U(u) \in \mathbb{R}$ is a convex function of $u$ and $F = F(u)$, naturally results from (\ref{eq:base}). The pair $(U, F)$ is often referred to as a generalized entropy-entropy flux pair. One such system is the one-dimensional Euler equations:
\begin{equation*}
    \frac{\partial}{\partial t} \begin{bmatrix}
    \rho \\ \rho u \\ \rho E
    \end{bmatrix} +
    \frac{\partial}{\partial x}
    \begin{bmatrix}
    \rho u \\ \rho u^2 + p \\ \rho u H
    \end{bmatrix} = 0,
\end{equation*}
where $\rho$ is the density, $u$ is the velocity, $E = e + \frac{1}{2}u^2$ is the total energy per unit mass, $\ p = \rho (\gamma - 1) e$ is the thermodynamic pressure, $H = E + \frac{p}{\rho} = \frac{a^2}{\gamma - 1} + \frac{1}{2}u^2$ is the total enthalpy and $a = \sqrt{\frac{\gamma p}{\rho}}$ is the speed of sound. An additional conservation equation that results from the above is that of entropy $\rho S$ ($S$ is the specific entropy):
\begin{equation*}
    \frac{\partial (\rho S)}{\partial t} + \frac{\partial (\rho u S)}{\partial x} = 0.
\end{equation*}
It is well known that discontinuous solutions to (\ref{eq:base}) can develop from smooth initial conditions. \textit{Weak} solutions must therefore be sought. Unfortunately, weak solutions are not uniquely defined and one needs additional conditions to distinguish physical numerical solutions from unphysical ones. For the Euler equations, solutions that introduce shock waves must lead to an increase in the entropy of the system. One way to enforce this property is to require the numerical solution to satisfy, in the sense of distributions:
\begin{equation}
    \frac{\partial (-\rho S)}{\partial t} + \frac{\partial (-\rho u S)}{\partial x} < 0.
\end{equation}
The negative sign on the left-hand side of the inequality is introduced as a convention to cast entropy production as a stability statement. The concept of a convex extension generalizes the concept of entropy to hyperbolic systems of conservation laws. It has been extensively studied and laid down by a number of researchers including Lax \cite{Lax, Lax1}, Friedrichs \cite{Lax}, Kruzkov \cite{Kruz}, Tadmor \cite{Tadmor3}, Harten \cite{Harten} and Godunov \cite{Godunov}. \\ 
%% PART 2: Tadmor's foundation and later extensions
\indent Numerical schemes that are consistent with the system (\ref{eq:base}) are not necessarily consistent with its convex extension (\ref{eq:base_entropy}). The pioneering work of Tadmor \cite{Tadmor5} introduced a class of finite-volume schemes that enforces either conservation or production of entropy at the semi-discrete level provided a certain condition is met by the interface flux. The first Entropy Conservative (EC) interface flux, which consists of a straight path integral in phase space, was also introduced. Entropy Stable (ES) schemes are achieved by combining an EC flux with a dissipation term that enforces production of entropy in the presence of discontinuities. Several developments followed. Of relevance to this work is the extension by Tadmor \cite{Tadmor1} and LeFloch \textit{et. al} \cite{LeFloch} of the framework to the fully discrete setting by studying the influence of time-integration on the entropy production. LeFloch \textit{et. al} \cite{LeFloch} introduced a Crank-Nicolson-like integration scheme using an intermediate state in time that introduces no production of entropy. One issue that hindered the use of EC/ES schemes in general is the non-closed form the first EC flux takes in the general nonlinear case (when the entropy is not a quadratic function of the state). Tadmor \cite{Tadmor1} introduced a new family of explicit EC fluxes which result from replacing the straight path along which the integral is computed into subpaths of an approximate Riemann solution in phase space. Roe \cite{Roe} introduced what is often referred to as an ``affordable'' entropy conservative flux for the Euler equations. This flux is obtained by solving the entropy conservation condition using algebraic manipulations inspired from \cite{Roe2}, and has been the baseline EC flux in many of the ES schemes that have been developed subsequently. \\ 
\indent The problem that motivates this work is the receding flow problem (also known as the \textit{123 problem} in Toro's book \cite{Toro}) studied extensively by Noh \cite{Noh} and Liou \cite{Liou, Liou2}. The problem consists of a smooth flow undergoing rarefaction caused by two flows receding from each other. It is identified by Liou as an open numerical problem as all well-known numerical fluxes give an anomalous temperature rise (often termed ``overheating'') at the origin that cannot be fixed by refining the mesh, decreasing the time-step, or increasing the solution order. Liou \cite{Liou2} first established a connection between the overheating and a spurious entropy rise after the first time step. In a more recent paper \cite{Liou}, he connected the entropy rise with the pressure component of the momentum flux. A cure proposed by Liou consists in replacing the energy equation with either the transport equation for the specific entropy $S$ or the conservation equation for the entropy $\rho S$. One issue with this approach is that conservation of total energy is no longer guaranteed. In view of this, the entropy conservative schemes introduced by Tadmor appear as an interesting option. These schemes, which enable the conservation of entropy without compromising the conservation of mass, momentum and energy, were not considered by Liou. An analysis similar to that of Liou \cite{Liou} within the framework of Tadmor is therefore worth carrying out. While the semi-discrete analysis of Liou suggests that the overheating could be avoided with an EC flux, a fully discrete analysis reveals the non-negligible influence of time-integration.  \\
\indent The effect of time-integration on the entropy conservation and entropy stability properties of a semi-discrete scheme has not been an active topic of research. From \cite{Tadmor1}, we know that if an entropy conservative flux is used in space, then the 1st-order explicit and implicit Euler schemes will make the resulting fully discrete system entropy unstable and stable, respectively. The question of the entropy stability of explicit and implicit schemes in general remains open. When the method of lines is used, Tadmor \cite{Tadmor1} shows that the entropy stability of the fully discrete scheme is determined by the balance between the entropy produced by the flux and the entropy produced by the time scheme. A fully discrete EC scheme is for instance obtained when an EC time scheme is used with an EC flux. To the authors knowledge, the only EC time scheme available to date is the one introduced by LeFloch et al. \cite{LeFloch}. However, this time scheme has issues very similar to the first EC flux (non-closed form) of Tadmor. This is no coincidence, as it turns out that the entropy conservation condition for spatial fluxes and the entropy conservation condition for this time scheme are similar. Building on this analogy and the algebraic manipulations that led to the EC Roe flux, we derive a new second-order EC time-integration scheme for the Euler equations and use it to show that conserving entropy at the fully discrete level does not necessarily prevent the overheating. \\   
 %% Paper organization
\indent The paper is organized as follows. In sections 2 and 3, we review Tadmor's numerical framework at the semi-discrete level. In section 4, we move to the fully discrete case: we examine the effect of time-integration on the entropy conservation properties of the underlying semi-discrete scheme, and introduce a new time-integration scheme for the Euler equations. In section 5, the overheating problem is analyzed in the context of EC schemes. We provide both analytical and numerical arguments explaining why the conservation of entropy is not exactly a remedy to the spurious behavior typically observed with finite-volume schemes.  

\section{The entropy variables}
An entropy-entropy flux pair $(U, F)$ is one that satisfies: 
\begin{equation}\label{eq:compa}
\frac{\partial U}{\partial u} \frac{\partial f}{\partial u} =\frac{\partial F}{\partial u} . 
\end{equation}
The entropy variables are defined by:
\begin{equation}
v = \bigg (\frac{\partial U}{\partial u}\bigg)^T.
\end{equation}
The convexity of the entropy function $U$ makes the mapping $u \rightarrow v$ one-to-one. So we can rewrite eq. (\ref{eq:base}) in terms of these variables:
\begin{equation*}
\frac{\partial u}{\partial t} + \frac{\partial g}{\partial x} = 0, \ g(v) = f(u(v)).
\end{equation*}
With the entropy variables, the system is symmetrized, meaning that the Jacobians of its temporal (denoted $H = H(v)$) and spatial fluxes are symmetric. Another property is that the temporal Jacobian is also positive definite. This implies that $u(v)$ and $g(v)$ can be written as gradients of potential functions $\phi$ and $\psi$:
\begin{equation*}
u(v) = \frac{\partial \phi}{\partial v}, \ g(v) = \frac{\partial \psi}{\partial v}.
\end{equation*}
With eq. (\ref{eq:compa}), it can be shown that these potential functions are given by:
\begin{equation*}
\phi(v) = v^Tu(v) - U(u(v)), \ \psi(v) = v^Tg(v) - G(v).
\end{equation*}

For the Euler equations, the entropy-entropy flux pair we work with is the one introduced by Hughes \textit{et. al} \cite{Hughes}:
\begin{equation*}
    U(u) = -\frac{\rho S}{\gamma - 1}, \ F(u) = -\frac{\rho u S}{\gamma - 1}, \ S = ln(p) - \gamma ln(\rho). 
\end{equation*}
This pair belongs to the class of entropy pairs derived by Harten \cite{Harten} and Tadmor \cite{Tadmor0}. This pair has interesting properties when applied to the Navier-Stokes equations: the Clausius-Duhem relation is recovered \cite{Hughes} and the corresponding potential functions are $ \phi = \rho, \psi = \rho u$. The corresponding entropy variables are given by:
\begin{equation*}
v = \bigg[\frac{\gamma - S}{\gamma - 1} - \frac{1}{2}\frac{\rho u^2}{p}, \ \frac{\rho u}{p}, \ -\frac{\rho}{p}\bigg]^T.
\end{equation*}
The temporal jacobian $H$ is given \cite{Barth, Hughes} by:
\begin{equation*}
    H = \begin{bmatrix}
    \rho   & \rho u       & \rho E \\
    \rho u & \rho u^2 + p & \rho H u \\
    \rho E & \rho H u     & \rho H^2 - \frac{a^2 p}{\gamma-1}
    \end{bmatrix}.
\end{equation*}

\section{The semi-discrete case}
\label{sec:SDEC}
Denote $j$ the cell index in space. A finite-volume discretization of eq. (\ref{eq:base}) is given by:
\begin{equation}\label{eq:FVM}
\frac{d}{dt}u_{j}(t) + \frac{1}{\Delta x}[f_{j + \frac{1}{2}} - f_{j - \frac{1}{2}}] = 0.
\end{equation}
In cell $j$, $u_{j}$ is the mean solution over the cell and $f_{j + \frac{1}{2}}$ is a consistent flux at the right interface. In the last section, we saw that for smooth (continuous) solutions of a system of conservation laws equipped with an entropy pair $(U, F)$, the entropy $U$ should be conserved according to eq. (\ref{eq:base_entropy}). For solutions with discontinuities, an admissibility criterion established by Lax \cite{Lax1} and Kruzkov \cite{Kruz} is that:
\begin{equation}\label{eq:entropy_stable}
\frac{\partial }{\partial t} U(u) + \frac{\partial}{ \partial x} F(u) \leq 0, 
\end{equation}
The semi-discrete scheme (\ref{eq:FVM}) is EC if it also solves a semi-discrete version of eq. (\ref{eq:base_entropy}):
\begin{equation}
\frac{d}{dt}U(u_{j}) + \frac{1}{\Delta x}[F_{j + \frac{1}{2}} - F_{j - \frac{1}{2}}] = 0,
\end{equation}
with $F_{j + \frac{1}{2}}$ a consistent entropy interface flux. Likewise, the scheme is entropy-stable if a semi-discrete version of eq. (\ref{eq:entropy_stable}):
\begin{equation}\label{eq:ESflux0}
\frac{d}{dt}U(u_{j}) + \frac{1}{\Delta x}[F_{j + \frac{1}{2}} - F_{j - \frac{1}{2}}] \leq 0,
\end{equation}
is inherently solved. Since
\begin{equation*}
v_{j}^T \frac{d}{dt} u_{j} = \frac{d}{dt} U(u_{j}),
\end{equation*}
an EC flux $f_{j + \frac{1}{2}}$ is one such that there exists a consistent entropy interface flux $F_{j + \frac{1}{2}}$ that satisfies:
\begin{equation*}
v_{j}^T [f_{j + \frac{1}{2}} - f_{j - \frac{1}{2}}] = F_{j + \frac{1}{2}} - F_{j - \frac{1}{2}}.
\end{equation*}
Tadmor \cite{Tadmor5} showed that this is equivalent to
\begin{equation}\label{eq:ECcond2}
[v_{j+1} - v_{j}]^T f_{j + \frac{1}{2}} = \psi_{j+1} - \psi_{j},
\end{equation}
where $\psi_{j} = \psi(u_j)$, and that if an interface flux $f_{j + \frac{1}{2}}$ is EC, then the corresponding interface entropy flux function $F_{j + \frac{1}{2}}$ is given by:
\begin{equation}\label{eq:entropy_interface}
F_{j + \frac{1}{2}} = \frac{1}{2}[v_{j} + v_{j+1}]^T f_{j + \frac{1}{2}} - \frac{1}{2}(\psi_{j} + \psi_{j + 1}).
\end{equation}

\subsection{Entropy Conservative fluxes}
The first EC flux was introduced by Tadmor \cite{Tadmor2}:
\begin{equation}\label{eq:TadmorEC}
f_{j + \frac{1}{2}} = \int_{0}^1 f(v_{j + \frac{1}{2}}(\xi)) d\xi, \  v_{j + \frac{1}{2}}(\xi) = v_{j} + \xi \Delta v_{j + \frac{1}{2}}, \ \Delta v_{j + \frac{1}{2}} = v_{j + 1} - v_{j}.
\end{equation}
This elegant flux has the inconvenient property of not having a closed form. Its evaluation requires using quadrature rules. The above flux is an integration along a straight path in phase space.
Later on, Tadmor \cite{Tadmor1} proposed to use piecewise-constant paths instead. The resulting flux has an explicit form which depends on the path decomposition. \\
%The neighboring grid values $v_{j}$ and $v_{j+1}$ are now connected through a set of intermediate states $\{ v_{j+\frac{1}{2}}^k\}_{k = 1}^N$. Let $\{ r_{j+\frac{1}{2}}^k\}_{k = 1}^N$ be a set of $N$ linearly independent vectors and let $\{ l_{j+\frac{1}{2}}^k\}_{k = 1}^N$ be the corresponding orthogonal vectors, i.e. the set that satisfies $\langle l_{j+\frac{1}{2}}^{k_1}, r_{j+\frac{1}{2}}^{k_2}\rangle = \delta_{k_1, k_2}$. This set is by definition the set of rows of the inverse of $R = \begin{bmatrix}  r_{j+\frac{1}{2}}^1 \ \dots \ r_{j+\frac{1}{2}}^N \end{bmatrix}$ . The piecewise constant path in phase space is then defined by:
%\begin{align*}
%\begin{cases}
%    v_{j + \frac{1}{2}}^1 =& \ v_{j}, \\
%    v_{j + \frac{1}{2}}^{k+1} =& \ v_{j+\frac{1}{2}}^k + \langle l_{j + \frac{1}{2}}^k, \Delta v_{j + \frac{1}{2}} \rangle r_{j+\frac{1}{2}}^k \\
%    v_{j + \frac{1}{2}}^N =& \ v_{j+1},
%\end{cases}
%\end{align*}
%and the corresponding flux is:
%\begin{equation*}
%f_{j + \frac{1}{2}} = \sum_{k = 1}^N \frac{\psi(v_{j + \frac{1}{2}}^{k+1}) - \psi(v_{j + \frac{1}{2}}^{k})}{\langle l_{j + \frac{1}{2}}^k, \Delta v_{j + \frac{1}{2}} \rangle} l_{j + \frac{1}{2}}^k.
%\end{equation*}
%Note that unlike the first version, this flux has an explicit form. The choice of the path is open.
\indent An EC flux that has been more popular over the years is the one derived by Roe \cite{Roe} for the Euler equations. The method used to derive it is central to this work. Denote $f^{*} = [f_1, \ f_2, \ f_3]$ the interface flux separating two adjacent cells. Using the jump notation, the condition (\ref{eq:ECcond2}) can be rewritten as:
\begin{equation*}
[v]^T f^* = [\psi].
\end{equation*}
For the Euler equations and the entropy pair we work with, the jump term in the entropy variables can be rewritten as:
\begin{equation*}
[v]^T = \bigg[\frac{ - [S]}{\gamma - 1} - \frac{1}{2}[\frac{\rho u^2}{p}], \ [\frac{\rho u}{p}], \ -[\frac{\rho}{p}]\bigg].
\end{equation*}
Define the set of independent variables $(z_1, z_2, z_3) = (\sqrt{\frac{\rho}{p}}, \sqrt{\frac{\rho}{p}}u, \sqrt{\rho p})$. Then
\begin{gather*}
    \rho = z_1 z_3, \ p = \frac{z_3}{z_1}, \ \frac{\rho}{p} = z_1^2, \ \frac{\rho u}{p} = z_1 z_2, \ \frac{\rho u^2}{p} = z_2^2, \ \rho u = z_2 z_3, \\
    S = (1 - \gamma) ln(z_3) - (1 + \gamma) ln(z_1).
\end{gather*}
Using the identities $ [a b] = \bar{a}[b] + \bar{b}[a]$ and $[ln(a)] = \frac{[a]}{a^{ln}}$, where $\bar{a}$ and $a^{ln}$ denote the arithmetic and logarithmic averages, respectively, one can show that:
\begin{gather*}
    [S] = \frac{(1 - \gamma)}{z_3^{ln}} [z_3] - \frac{(1 + \gamma)}{z_1^{ln}} [z_1], \
    [\frac{\rho u^2}{p}] = 2\bar{z_2}[z_2], \
    [\frac{\rho u}{p}] = \bar{z_1}[z_2] + \bar{z_2}[z_1], \\
    [\frac{\rho}{p}] = 2\bar{z_1}[z_1], \ [\rho u] = \bar{z_2}[z_3] + \bar{z_3}[z_2].
\end{gather*}
The motivation behind the introduction of the variables $z_i$ is to ``exactly linearize'' all the jump terms involved in the entropy conservation equation (\ref{eq:ECcond2}). The entropy conservation condition  now becomes:
\begin{equation*}
    f_1 ( \frac{1}{z_3^{ln}} [z_3] - \frac{1 + \gamma}{1 - \gamma}\frac{1}{z_1^{ln}} [z_1] - \bar{z_2}[z_2]) +
    f_2 (\bar{z_1}[z_2] + \bar{z_2}[z_1]) + f_3 (- 2 \bar{z_1}[z_1]) = \bar{z_2}[z_3] + \bar{z_3}[z_2].
\end{equation*}
Regrouping, the above becomes:
\begin{equation*}
    [z_1](- f_1 \frac{1 + \gamma}{1 - \gamma}\frac{1}{z_1^{ln}} + f_2 \bar{z_2} - 2 f_3 \bar{z_1}) + [z_2] ( - f_1\bar{z_2} + f_2 \bar{z_1}) + [z_3](\frac{1}{z_3^{ln}} f_1)
     = [z_2] \bar{z_3} + [z_3]\bar{z_2}.
\end{equation*}
The jumps in the $z_i$ are independent, therefore:
\begin{align*}
- f_1 \frac{1 + \gamma}{1 - \gamma}\frac{1}{z_1^{ln}} + f_2 \bar{z_2} - 2 f_3 \bar{z_1} =& \ 0, \
 - f_1\bar{z_2} + f_2 \bar{z_1} = \bar{z_3}, \
 \frac{1}{z_3^{ln}} f_1 = \bar{z_2}. 
\end{align*}
The variables $z_i$ basically enabled us to convert the scalar condition (\ref{eq:ECcond2}) into a system of 3 equations that we can easily solve:
\begin{align*}
f_1 =& \ \bar{z_2} z_3^{ln} , \
f_2 = (\bar{z_3} + f_1 \bar{z_2})/(\bar{z_1}) , \
f_3 = \frac{1}{2 \bar{z_1}}(- f_1 \frac{1 + \gamma}{1 - \gamma}\frac{1}{z_1^{ln}} + f_2 \bar{z_2}).
\end{align*}
This is the Roe EC flux. \\
\indent The choice of the variables $z_i$ is open. Using the same method with the set $(z_1, z_2, z_3) = (\rho, u, \frac{\rho}{2p})$, Chandrasekhar \cite{Chandra} derived the following EC flux: 
\begin{align*}
    f_1 =& \ z_1^{ln}\bar{z_2}, \ f_2 = \frac{\bar{z_1}}{2\bar{z_3}} + \bar{z_2}f_1, \    f_3 = \bigg[\frac{1}{2(\gamma - 1)z_3^{ln}} - \frac{1}{2}\bar{z_2^2}\bigg]f_1 + \bar{z_2}f_2.
\end{align*}
Jameson \cite{Jameson} showed that an interface flux preserves the kinetic energy of the system at the semi-discrete level provided that the density flux $f_1$ and momentum flux $f_2$ satisfy $f_2 = \tilde{p} + \overline{u} f_1$, where $\tilde{p}$ is any consistent average pressure. The above EC flux satisfies this property as well. In contrast to the conclusions of \cite{Chandra} (section 4.6.), such a flux is not unique. With the set $ (z_1, z_2, z_3) = (p, u, \frac{\rho}{2p})$, the resulting EC flux:
\begin{align*}
f_1 =& \ 2 \bar{z_3}\bar{z_2} z_1^{ln}, \
f_2 = \frac{\bar{\rho}}{2 \bar{z_3}} + \bar{u} f_1, \
f_3 = f_1 (\frac{\gamma}{\gamma - 1}\frac{1}{2 z_3^{ln}} - \frac{1}{2}\overline{z_2^2}) + f_2 \bar{z_2} - \bar{z_1} \bar{z_2},
\end{align*}
is also kinetic energy preserving. The term $-\bar{z_1} \bar{z_2}$ of the energy flux $f_3$ is missing in \cite{Chandra}. \\ 
\indent Although the EC Roe flux was introduced for the Euler equations, the methodology used to derive it is general and can be applied to other systems in multiple dimensions such as the ideal MagnetoHydroDynamics (MHD) equations (see Winters \textit{et al.} \cite{Winters}).
\subsection{Entropy stable spatial fluxes}
Entropy stable fluxes are fluxes such that inequality (\ref{eq:ESflux0}) is enforced, with $F_{j + \frac{1}{2}}$ a consistent entropy interface flux. Entropy Stable (ES) fluxes are typically built by combining an entropy stabilization term $f_{j + \frac{1}{2}}^S$ with an EC flux $f_{j + \frac{1}{2}}^*$. An ES flux has the general form $f_{j + \frac{1}{2}} = f_{j + \frac{1}{2}}^* - f_{j + \frac{1}{2}}^S$ and the entropy stabilization term should be such that:
\begin{equation}\label{eq:ESflux}
 v_{j}^T [f_{j + \frac{1}{2}} - f_{j - \frac{1}{2}}] = [F_{j + \frac{1}{2}} - F_{j - \frac{1}{2}}] + \mathcal{E}, \ \ \mathcal{E} > 0,
\end{equation}
so that, at the semi-discrete level:
\begin{equation*}
\frac{d}{dt}U(u_{j}(t)) + \frac{1}{\Delta x}[F_{j + \frac{1}{2}} - F_{j - \frac{1}{2}}] = -\mathcal{E} < 0.
\end{equation*}
Note that the entropy flux $F_{j + \frac{1}{2}}$ does not necessarily need to be the one generated by $f_{j + \frac{1}{2}}^*$  (denoted $F_{j + \frac{1}{2}}^*$). Tadmor \cite{Tadmor2} showed that if the stabilizer has the general form
\begin{equation*}
f_{j + \frac{1}{2}}^S = Q_{j + \frac{1}{2}} \Delta v_{j + \frac{1}{2}},
\end{equation*}
with $Q_{j + \frac{1}{2}} \in \mathbb{R}^{N \times N}$, then the condition (\ref{eq:ESflux}) is almost met with:
\begin{equation*}
F_{j + \frac{1}{2}} = F_{j + \frac{1}{2}}^* - v_{j + 1}^T Q_{j + \frac{1}{2}} \Delta v_{j + \frac{1}{2}},
\end{equation*}
and 
\begin{equation}\label{eq:SpaceDiss}
\mathcal{E} = \frac{1}{\Delta x}\big[ \Delta v_{j + \frac{1}{2}}^T Q_{j + \frac{1}{2}} \Delta v_{j + \frac{1}{2}} + \Delta v_{j - \frac{1}{2}}^T Q_{j - \frac{1}{2}} \Delta v_{j - \frac{1}{2}} \big].
\end{equation}
Entropy stability is achieved if the dissipation matrix $Q_{j + \frac{1}{2}}$ positive definite. 
%Note that eq. (\ref{eq:SpaceDiss}) tells us exactly how much entropy is dissipated in each element as a consequence of the ES flux that is used in space.    
\section{The fully discrete case}
\subsection{Entropy stability of time schemes}
Let's assume that an EC flux is used in space. What happens at the fully discrete level, when time is discretized? If we evolve in time using Forward Euler (FE) for instance:
\begin{equation*}
u_{j}^{n+1} - u_{j}^{n} + \lambda [f_{j + \frac{1}{2}}^{n} - f_{j - \frac{1}{2}}^{n}] = 0, \ \lambda = \frac{\Delta t}{\Delta x},
\end{equation*}
are we simultaneously evolving in time the entropy equation in with FE as well? In other words, is the above fully discrete scheme also solving
\begin{equation*}
U(u_{j}^{n+1}) - U(u_{j}^{n}) + \lambda [F_{j + \frac{1}{2}}^{n} - F_{j - \frac{1}{2}}^{n}] = 0 \  \mbox{?}
\end{equation*}
In the above and what follows, the superscript $n$ refers to the time instant. We know that
\begin{equation*}
    (v_{j}^n)^T [f^n_{j + \frac{1}{2}} - f^n_{j - \frac{1}{2}}] = F^n_{j + \frac{1}{2}} - F^n_{j - \frac{1}{2}},
\end{equation*}
therefore the answer depends on whether
\begin{equation*}
(v_{j}^{n})^T [u_{j}^{n+1} - u_{j}^{n}] = U(u_{j}^{n+1}) - U(u_{j}^{n})
\end{equation*}
holds. Tadmor \cite{Tadmor1} showed that, for FE:
\begin{subequations}
\begin{align}
(v_{j}^{n})^T [u_{j}^{n+1} - u_{j}^{n}] = U(u_{j}^{n+1}) - U(u_{j}^{n}) - \mathcal{E}^{FE}(v_{j}^{n}, v_{j}^{n+1}), \label{eq:TadFE}\\
\mathcal{E}^{FE}(v_{j}^{n}, v_{j}^{n+1}) = \int_{-\frac{1}{2}}^{\frac{1}{2}} (\frac{1}{2} + \xi)(\Delta v_{j}^{n + \frac{1}{2}})^T H(v_{j}^{n + \frac{1}{2}}(\xi)) \Delta v_{j}^{n + \frac{1}{2}} d\xi > 0,
\end{align}
\end{subequations}
where $v_{j}^{n + \frac{1}{2}}(\xi) = \frac{v_{j}^{n+1} + v_{j}^{n}}{2} + \xi (v_{j}^{n+1} - v_{j}^{n})$ and $
\Delta v_{j}^{n + \frac{1}{2}} = v_{j}^{n+1} - v_{j}^{n}$. This means that at the fully discrete level:
\begin{equation}\label{eq:FDFE}
U(u_{j}^{n+1}) - U(u_{j}^{n}) + \lambda [F_{j + \frac{1}{2}}^{n} - F_{j + \frac{1}{2}}^{n}] = \mathcal{E}^{FE} > 0.
\end{equation} 
This makes FE unconditionally entropy unstable. For Backward Euler, Tadmor \cite{Tadmor1} showed that:
\begin{subequations}
\begin{align}
(v_{j}^{n+1})^T [u_{j}^{n+1} - u_{j}^{n}] = U(u_{j}^{n+1}) - U(u_{j}^{n}) + \mathcal{E}^{BE}(v_{j}^{n}, v_{j}^{n+1}), \label{eq:TadBE}\\
\mathcal{E}^{BE}(v_{j}^{n}, v_{j}^{n+1}) = \int_{-\frac{1}{2}}^{\frac{1}{2}} (\frac{1}{2} - \xi)(\Delta v_{j}^{n + \frac{1}{2}})^T H(v_{j}^{n + \frac{1}{2}}(\xi)) \Delta v_{j}^{n + \frac{1}{2}} d\xi > 0.
\end{align}
\end{subequations}
This means that at the fully discrete level:
\begin{equation}\label{eq:FDBE}
U(u_{j}^{n+1}) - U(u_{j}^{n}) + \lambda [F_{j + \frac{1}{2}}^{n+1} - F_{j + \frac{1}{2}}^{n+1}] = - \mathcal{E}^{BE} < 0
\end{equation}
This makes BE unconditionally entropy stable. One may wonder if all implicit and explicit time-integration schemes are unconditionally entropy stable and unstable, respectively. To the best of the authors' knowledge, this is an open question. To support this statement, let's use the two main results of Tadmor's analysis, i.e. eqns. (\ref{eq:TadFE}) and (\ref{eq:TadBE}),
to derive the entropy production of some well-known time-integration schemes. Define:
\begin{equation*}
    R_j^f(u) = -\frac{1}{\Delta x} (f_{j+\frac{1}{2}} - f_{j-\frac{1}{2}}), \ R_j^F(u) = v_j^T R_j^f(u) = -\frac{1}{\Delta x} (F_{j+\frac{1}{2}} - F_{j-\frac{1}{2}}),
\end{equation*}
so that we can rewrite eq. (\ref{eq:FVM}) in a more compact form. 
The implicit 2nd-order backward difference scheme is given by:
\begin{equation*}
    u_{j}^{n+2} - \frac{4}{3}u_{j}^{n+1}+\frac{1}{3}u_{j}^n = \frac{2}{3}\Delta t R_j^f(u^{n+2}).
\end{equation*}
If we rewrite the above scheme as:
\begin{equation*}
    \frac{4}{3} (u_{j}^{n+2} - u_{j}^{n+1}) - \frac{1}{3}(u_{j}^{n+2} - u_{j}^n) = \frac{2}{3}\Delta t R_j^f(u^{n+2}),
\end{equation*}
left-multiply by $(v_j^{n+2})^T$ and use eq. (\ref{eq:TadBE}), we obtain the following for the discrete entropy:
\begin{align*}
    U(u_{j}^{n+2}) - \frac{4}{3}U(u_{j}^{n+2}) + \frac{1}{3} U(u_{j}^n) = \frac{2}{3}\Delta t R_j^F(u^{n+2}) - \mathcal{E}^{BDF2}
\end{align*}
This is basically BDF2 for the discrete entropy, with an additional entropy production term $\mathcal{E}^{BDF2}$ given by:
\begin{equation*}
    \mathcal{E}^{BDF2} = \frac{4}{3} \mathcal{E}^{BE}(v_{j}^{n+1},v_{j}^{n+2}) - \frac{1}{3} \mathcal{E}^{BE}(v_{j}^{n},v_{j}^{n+2})    
\end{equation*}
The production term $\mathcal{E}^{BE}(v_{j}^{n},v_{j}^{n+2})$ can determine the entropy stability of BDF2. However, its sign is hard to establish. The explicit Leap-Frog Method is given by:
\begin{equation*}
u_{j}^{n+1} = u_{j}^{n-1} + 2 \Delta t R_j^f(u^n).
\end{equation*}
Subtracting $u_{j}^n$ on both sides, left-multiplying by $(v_{j}^n)^T$, and using eqns. (\ref{eq:TadFE}) and (\ref{eq:TadBE}), we get a Leap-Frog of the entropy
\begin{align*}
U(u_{j}^{n+1}) = U(u_{j}^{n-1}) + 2 \Delta t R_j^F(u^n) - \mathcal{E}^{LF}
\end{align*}
with an entropy production term $\mathcal{E}^{LF}$ given by:
\begin{equation*}
    \mathcal{E}^{LF} = \mathcal{E}^{BE}(v_{j}^{n-1}, v_{j}^{n}) - \mathcal{E}^{FE}(v_{j}^n, v_{j}^{n+1}).
\end{equation*}
Here again, it is hard to make a statement about the sign of the entropy production term $\mathcal{E}^{LF}$. We could derive similar results for other schemes and reach the same conclusion. 

\subsection{Entropy conservative time schemes}
We still assume that an EC flux is used in space. LeFloch \textit{et al.} \cite{LeFloch} (Theorem 3.1.) showed that the following scheme:
\begin{equation}\label{eq:LeFloch}
    \frac{u_j^{n+1} - u_j^{n}}{\Delta t} = R_j^f(u(v^{n + \frac{1}{2}})),
\end{equation}
with $v^{n + \frac{1}{2}}$ an intermediate state in time given by: 
\begin{equation}\label{eq:TECTadmor}
v^{n + \frac{1}{2}} = \int_{-\frac{1}{2}}^{\frac{1}{2}} v(\frac{u^{n} + u^{n+1}}{2} + \xi \Delta u^{n + \frac{1}{2}}) d\xi, \Delta u^{n + \frac{1}{2}} = u^{n+1} - u^{n},
\end{equation}
is entropy conservative in the sense that the scheme satisfies: 
\begin{equation*}
    \frac{U(u_j^{n+1}) - U(u_j^{n})}{\Delta t} = R_j^F(u(v^{n + \frac{1}{2}})).
\end{equation*}
The scheme is also shown to be second-order accurate, in the sense that eq. (\ref{eq:LeFloch}) is a second-order approximation to eq. (\ref{eq:base}) evaluated at $t = \frac{t^{n+1} + t^{n}}{2}$.
If the system of interest is symmetric ($v = u$), the intermediate state becomes $ v^{n + \frac{1}{2}} = \frac{u^{n} + u^{n+1}}{2}$.

Tadmor refers to this scheme as a ``Generalized Crank-Nicolson'' scheme in \cite{Tadmor1}. In the general nonlinear case, this scheme is ``inpractical'' to the same extent as the first EC flux. The intermediate state does not have a closed form and requires quadrature. \\
\indent The similarity between the intermediate state in time (\ref{eq:TECTadmor}) and the first EC flux (\ref{eq:TadmorEC}) is no coincidence. The condition on the intermediate state for the proposed scheme to be entropy conservative is (see Assumption 2.1. in \cite{LeFloch} for $q = 1$):
\begin{equation}\label{eq:TIEC}
(v_{j}^{n+\frac{1}{2}})^T [u_{j}^{n+1} - u_{j}^{n}] = U(u_{j}^{n+1}) - U(u_{j}^{n}).
\end{equation}
This equation is an analog of the Entropy conservation condition in space (\ref{eq:ECcond2}). We can therefore apply the technique used to derive the EC Roe flux to derive an ``affordable'' intermediate state in time. Denote $v^{n + \frac{1}{2}} = [v_1, \ v_2, \ v_3]$. 
%The jumps in time are in terms of the quantities:  
%\begin{align*}
 %   u = [\rho, \ \rho u, \rho E]^T, \ U = \frac{-\rho S}{\gamma - 1}. 
%\end{align*}
Let's consider $\rho, \ u$ and $p$ as the independent variables. The jumps can be written as:
\begin{gather*}
    [\rho u] = \bar{\rho} [u] + \bar{u} [\rho], \
    [\rho E] = \frac{[p]}{\gamma - 1} + \frac{1}{2} [\rho u^2] = \frac{[p]}{\gamma - 1} + \frac{1}{2} \bar{u^2} [\rho] + \bar{\rho} \bar{u} [u], \\
    [\rho S] = \bar{\rho} [S] + \bar{S} [\rho] 
             = \bar{\rho} \frac{[p]}{p^{ln}} - \gamma \bar{\rho} \frac{[\rho]}{\rho^{ln}} + \bar{S} [\rho].
\end{gather*}
Injecting the above in eq. (\ref{eq:TIEC}) and regrouping we obtain:
\begin{equation*}
    [\rho](v_1 + \bar{u}v_2 + \frac{1}{2}\bar{u^2}v_3) + [u](\bar{\rho}v_2 + \bar{u}\bar{\rho} v_3) + [p](\frac{v_3}{\gamma - 1}) = \frac{-1}{\gamma - 1}( [\rho] (\bar{S} - \gamma \frac{\bar{\rho}}{\rho^{ln}}) + [p]\frac{\bar{\rho}}{p^{ln}} ). 
\end{equation*}
The jumps are independent, therefore:
\begin{align*}
    v_1 + \bar{u}v_2 + \frac{1}{2}\bar{u^2}v_3 =& \ \frac{1}{\gamma -1}( \gamma - \bar{S}  \frac{\bar{\rho}}{\rho^{ln}}), \
    \bar{\rho}v_2 + \bar{u}\bar{\rho} v_3 = 0, \
    v_3 = - \frac{\bar{\rho}}{p^{ln}}. 
\end{align*}
The solution is:
\begin{align*}
    v_1 =& \ \frac{1}{\gamma - 1}( \gamma \frac{\bar{\rho}}{\rho^{ln}} -\bar{S} ) -\bar{u}v_2 - \frac{1}{2}\bar{u^2} v_3 , \
    v_2 = - \bar{u}v_3, \
    v_3 = - \frac{\bar{\rho}}{p^{ln}}.
\end{align*}
This intermediate state satisfies condition (\ref{eq:TIEC}) and is consistent. Let's show that the resulting scheme is second-order as well. A Taylor expansion about $t = \frac{t^{n} + t^{n+1}}{2}$ gives:
\begin{equation*}
    \frac{u_j^{n+1} - u_j^{n}}{\Delta t} = \partial_t u (x_j, \frac{t^{n} + t^{n+1}}{2}) + \mathcal{O}(\Delta t^2).
\end{equation*}
To conclude, let's show that:
\begin{equation*}
    v^{n+\frac{1}{2}} = v(\frac{t^n + t^{n+1}}{2}) + \mathcal{O}(\Delta t^2).
\end{equation*}
Let's establish a few results first. Let $a(t)$ be a strictly positive time-dependent quantity and denote $a^n = a(t^n)$, $a^{n+1} = a(t^{n+1})$ and $a^* = a(\frac{t^n + t^{n+1}}{2})$. Using a Taylor analysis about $t = \frac{t^n + t^{n+1}}{2} = t^n + \frac{\Delta t}{2} = t^{n+1} - \frac{\Delta t}{2}$, one can show that:
\begin{gather*}
    \bar{a} = a^* + \mathcal{O}(\Delta t^2), \
    \bar{a^2} = (a^*)^2 + \mathcal{O}(\Delta t^2), \
    a^{n+1} - a^n = \Delta t a^*_{,t} + \mathcal{O}(\Delta t^3), \\ log(a^{n+1}) - log(a^{n}) = \Delta t \frac{a^*_{,t}}{a^*} + \mathcal{O}(\Delta t^3).
\end{gather*}
Therefore:
\begin{gather*}
    a^{ln} = \frac{a^{n+1} - a^n}{log(a^{n+1}) - log(a^{n})} = \frac{\Delta t a^*_{,t} + \mathcal{O}(\Delta t^3)}{\Delta t \frac{a^*_{,t}}{a^*} + \mathcal{O}(\Delta t^3)} =  \frac{a^* + \mathcal{O}(\Delta t^2)}{1 + \mathcal{O}(\Delta t^2)} 
    = a^* + \mathcal{O}(\Delta t^2)
\end{gather*}
%We used the result
%\begin{equation*}
%    \frac{1}{1 + \mathcal{O}(\Delta t^2)} = 1 + %\mathcal{O}(\Delta t^2),
%\end{equation*}
%that is a consequence of:
%\begin{equation*}
%    \frac{1}{1 + x} = 1 - x + \mathcal{O}(x^2), \ x \rightarrow 0.
%\end{equation*}
Likewise we show another useful identity:
\begin{equation*}
    \frac{\bar{a}}{b^{ln}} = \frac{a^*}{b^*} + \mathcal{O}(\Delta t^2),
\end{equation*}
where $b(t)$ is another strictly positive quantity. With all the above results we can show that our nonlinear intermediate state $v_{n+\frac{1}{2}} = [v_1, \ v_2, \ v_3]$ satisfies:
\begin{gather*}
    v_1 = \frac{\gamma - S^*}{\gamma - 1} -\frac{1}{2}\frac{\rho^* (u^*)^2}{p^*} + \mathcal{O}(\Delta t^2), \
    v_2 = \frac{\rho^* u^*}{p^*} + \mathcal{O}(\Delta t^2), \
    v_3 = - \frac{\rho^*}{p^*} + \mathcal{O}(\Delta t^2).
\end{gather*}
Note that this methodology can be applied to other hyperbolic systems with a convex extension (shallow water equations, ideal MHD equations, etc...).
\section{The overheating problem: receding flow}
\begin{figure}[h!]
    \centering
    \subfigure[Density]{\includegraphics[scale =
    0.60]{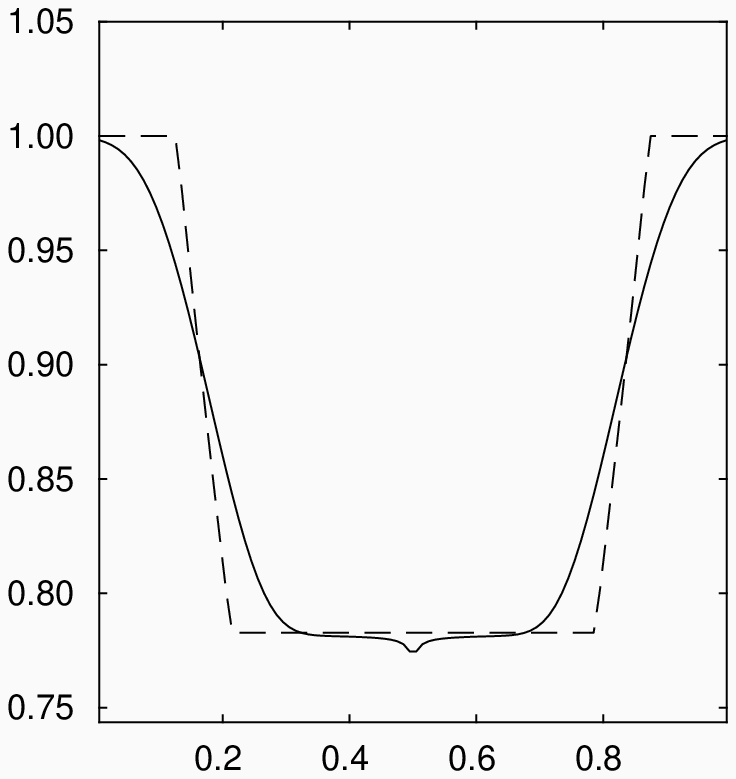}}
    \subfigure[Pressure]{\includegraphics[scale = 0.60]{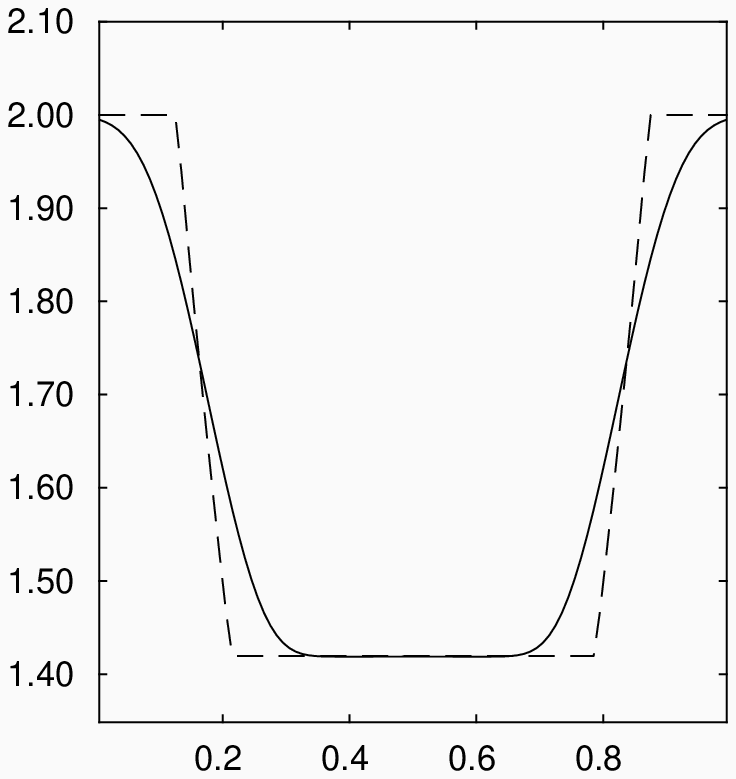}}
    \subfigure[Velocity]{\includegraphics[scale = 0.60]{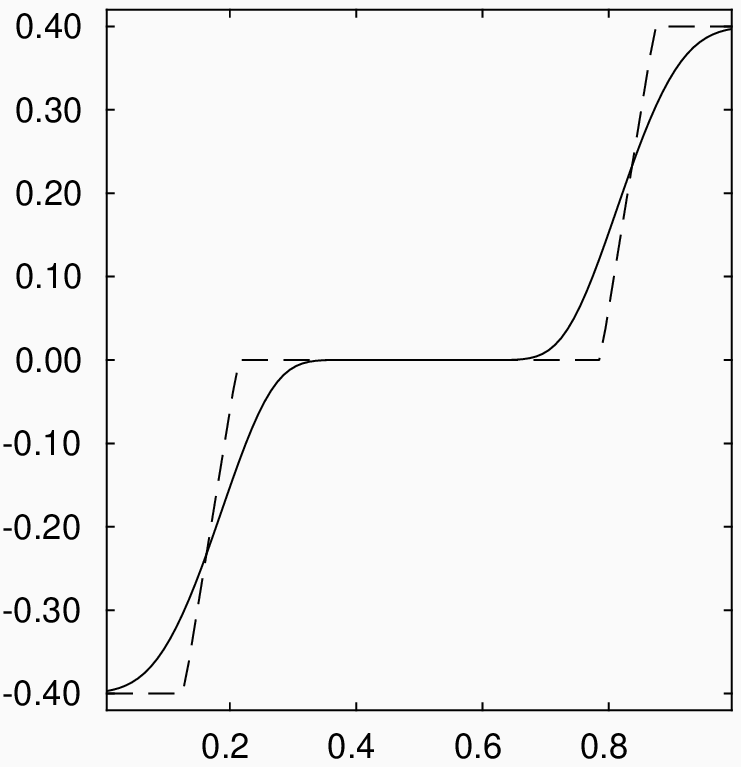}}
    \subfigure[Specific Entropy]{\includegraphics[scale = 0.60]{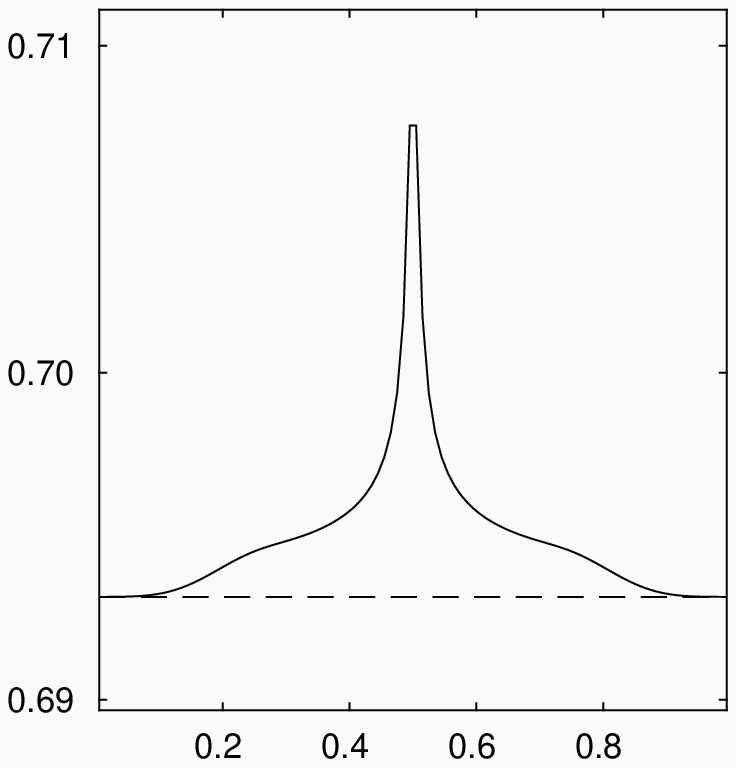}}
    \caption{Receding flow problem: Numerical solution (full line) at $t = 0.18s$ with the Roe flux and Forward Euler in time. 100 cells and $\Delta t = 10^{-3}s$. }
    \label{fig:RF_classicRoe}
\end{figure}

The receding flow problem \cite{Liou, Liou2} is a 1D Riemann problem defined by the following initial conditions:
\begin{equation}
    u_L < 0, \ u_R > 0, \rho_L = \rho_R = \rho^0, \ p_L = p_R = p^0.
\end{equation}
where the subscripts L and R refer to the left and right sides of the domain, respectively. Liou considered the case of equal velocity magnitudes $|u_L| = |u_R| = u^0$. \\
\indent Liou describes this problem as ``fundamental'' in the sense that the overheating cannot be overcome by refining the mesh or changing the time step (it is independent of the CFL number). One of the main findings of his study is that the overheating originates from an \textit{ab initio} entropy production at the beginning of the run. Figure \ref{fig:RF_classicRoe} (Roe flux in space, forward euler in time) illustrates the numerical behavior that is typically observed with the wide range of fluxes Liou \cite{Liou, Liou2} considered. The pressure is well resolved whereas the density is slightly under-estimated at the center (see figure \ref{fig:RF_classicRoe}-(a)). That the overheating is generated at the very first instant is intuitive given the nature of rarefaction waves (discontinuities that should vanish after some time) in contrast to shock waves (discontinuities that persist in time).

\subsection{Liou's semi-discrete analysis}
\begin{figure}[h!]
    \centering
    \subfigure[Density]{\includegraphics[scale = 0.63]{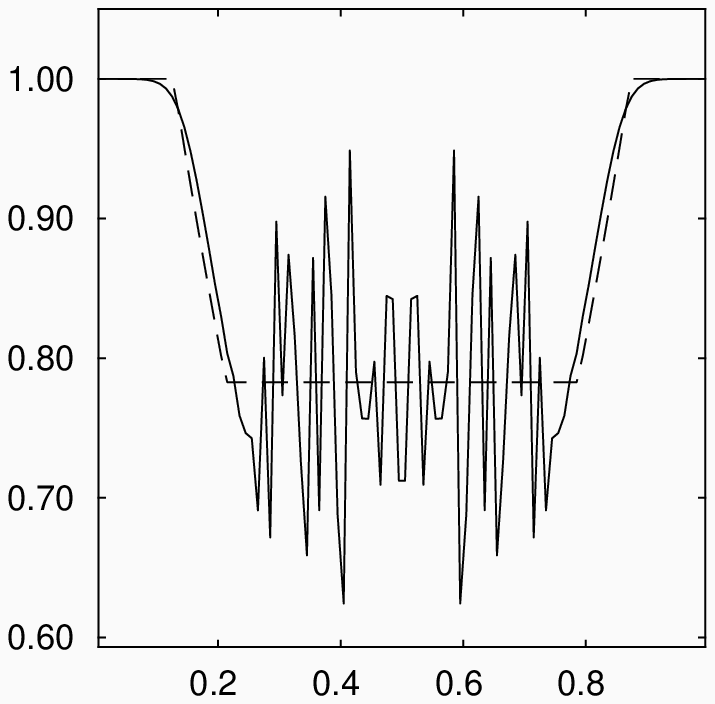}}
    \subfigure[Pressure]{\includegraphics[scale = 0.63]{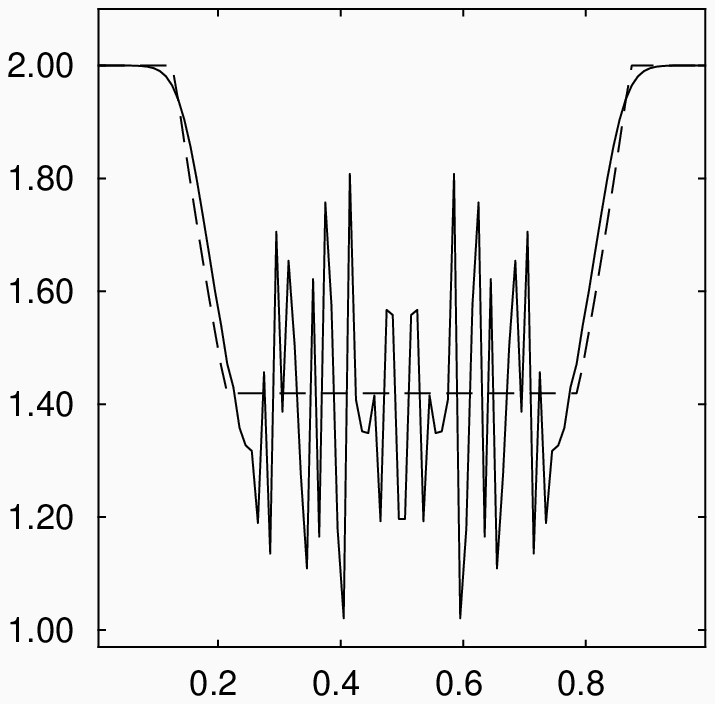}}
    \subfigure[Velocity]{\includegraphics[scale = 0.63]{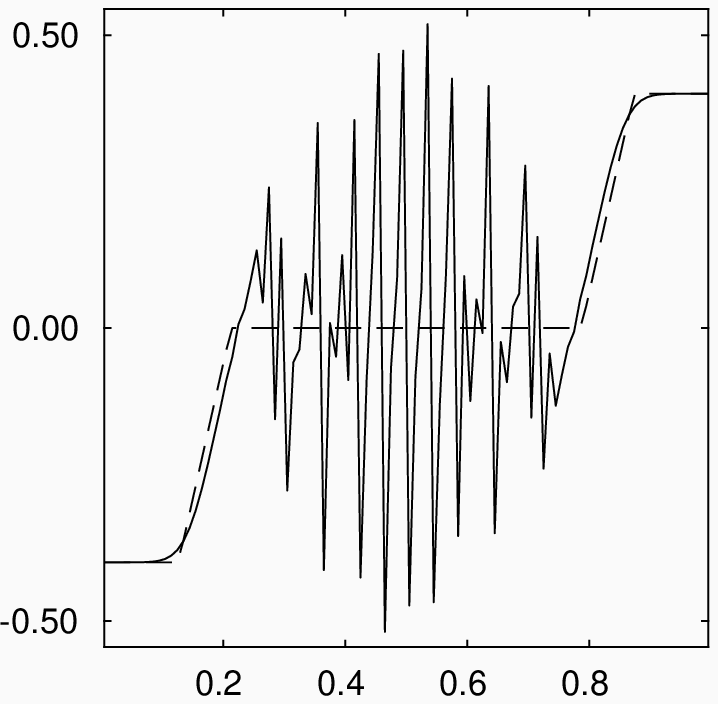}}
    \subfigure[Specific Entropy]{\includegraphics[scale = 0.63]{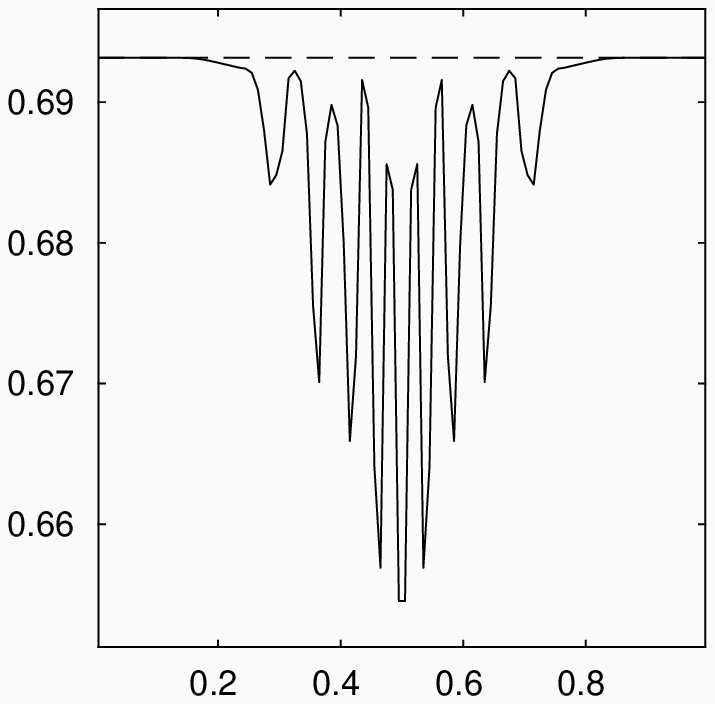}}
    \caption{Receding flow problem: Numerical solution (full line) at $t = 0.18s$ with the EC Roe flux and Forward Euler in time. 100 cells and $\Delta t = 10^{-3}s$. }
    \label{fig:RF_FEECRoe}
\end{figure}

To investigate how $S$ is produced in the discretized conservation laws, Liou \cite{Liou} begins with the following equation:
\begin{equation*}
    \frac{p}{R}\frac{\partial S}{\partial t} =  - (\frac{a^2}{\gamma - 1} - \frac{u^2}{2})\frac{\partial \rho}{\partial t} - u \frac{\partial \rho u}{\partial t}  + \frac{\partial \rho E}{\partial t}.
\end{equation*}
R is the gas constant. This equation relates the temporal variation of the specific entropy $S$ to that of mass, momentum and total energy. \\
\indent Denote cell ``R'' as the cell immediately to the right of the interface with index 1. Then integration over cell R gives:
\begin{multline*}
    \oint_R \frac{p}{R}\frac{\partial S}{\partial t}dV =   (\frac{(a^0)^2}{\gamma - 1} - \frac{(u^0)^2}{2})[(\rho u)_{3/2} - (\rho u)_{1/2}] + u^0 [(\rho u^2 + p)_{3/2} - (\rho u^2 + p)_{1/2}] \\ - [(\rho u H)_{3/2} - (\rho u H)_{1/2}].
\end{multline*}
The fluxes at interface $3/2$ are determined by the initial conditions:
\begin{equation*}
    (\rho u)_{3/2} = \rho^0 u^0, \ (\rho u^2 + p)_{3/2} = \rho^0 (u^0)^2 + p^0, \ (\rho u H)_{3/2} = \rho^0 u^0 H^0.
\end{equation*}
For all the fluxes tested by Liou \cite{Liou}, the values at the interface $1/2$ are given by:
\begin{equation*} 
    (\rho u)_{1/2} = 0, \ (\rho u^2 + p)_{1/2} = m_{1/2} + p_{1/2}, \ (\rho u H)_{1/2} = 0.
\end{equation*}
$m_{1/2}$ and $p_{1/2}$ are the momentum and pressure fluxes, respectively. Combining the above 3 equations results in the following:
\begin{equation*}
    \oint_R \frac{p}{R}\frac{\partial S}{\partial t}dV = u_0 [(p^0 - p_{1/2}) - m_{1/2}].
\end{equation*}
The right-hand side term remains non-zero for all the fluxes Liou considered. An identical result is obtained for the ``L'' cell immediately to the left of the interface, meaning that the entropy rise occurs symmetrically about the interface. In light of the above equation, Liou attributed the \textit{ab initio} generation of entropy to the pressure and momentum components of the numerical flux. He concluded his study by showing that replacing the energy equation with the conservation equation for entropy cures the overheating. One of the issues with this cure is the lack of conservation of energy. \\
\indent Liou's study did not consider Tadmor's family of schemes. EC schemes allow the conservation of entropy without compromising that of energy. If the EC Roe flux is used, the flux values at interface $1/2$ take the values:
\begin{equation*}
    (\rho u)_{1/2} = 0, \ (\rho u^2 + p)_{1/2} = p^0, \ (\rho u H)_{1/2} = 0.
\end{equation*}
and we obtain, for both the R and L cells:
\begin{equation*}
    \oint \frac{p}{R}\frac{\partial S}{\partial t}dV = 0.
\end{equation*}
This suggests that the spurious entropy production would be avoided. This is unfortunately not the case (see figure \ref{fig:RF_FEECRoe}). In the next section, we use a fully discrete analysis to explain why the spurious entropy rise is not necessarily avoided.
\begin{figure}[h!]
    \centering
    \subfigure[Density]{\includegraphics[scale = 0.63]{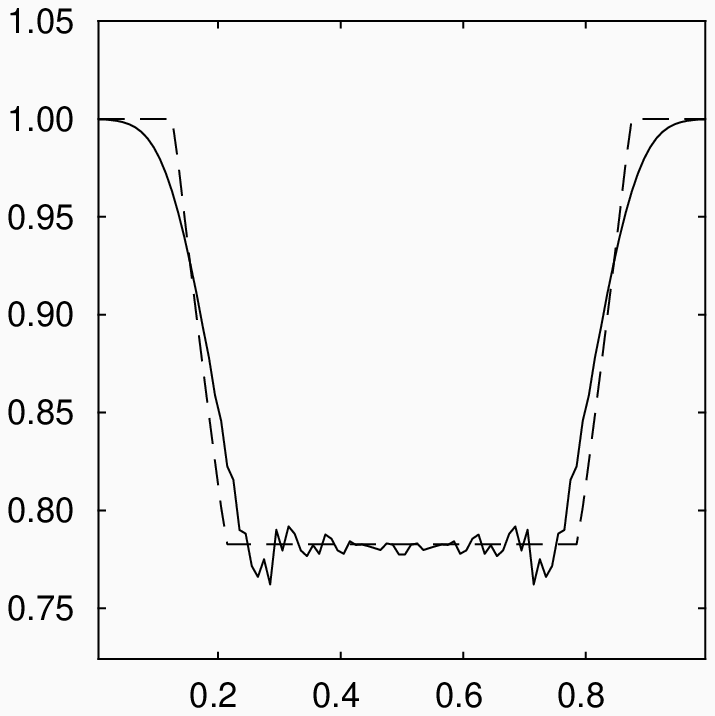}}
    \subfigure[Pressure]{\includegraphics[scale = 0.63]{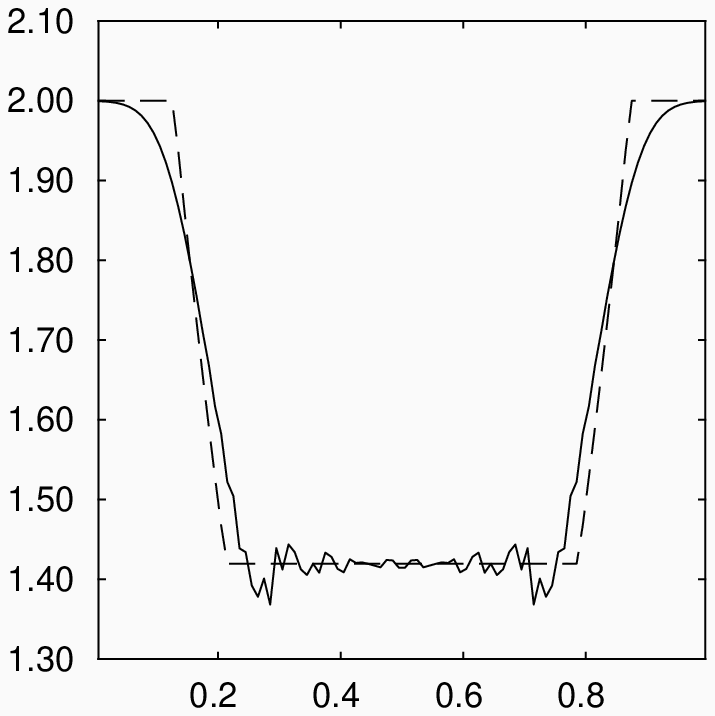}}
    \subfigure[Velocity]{\includegraphics[scale = 0.63]{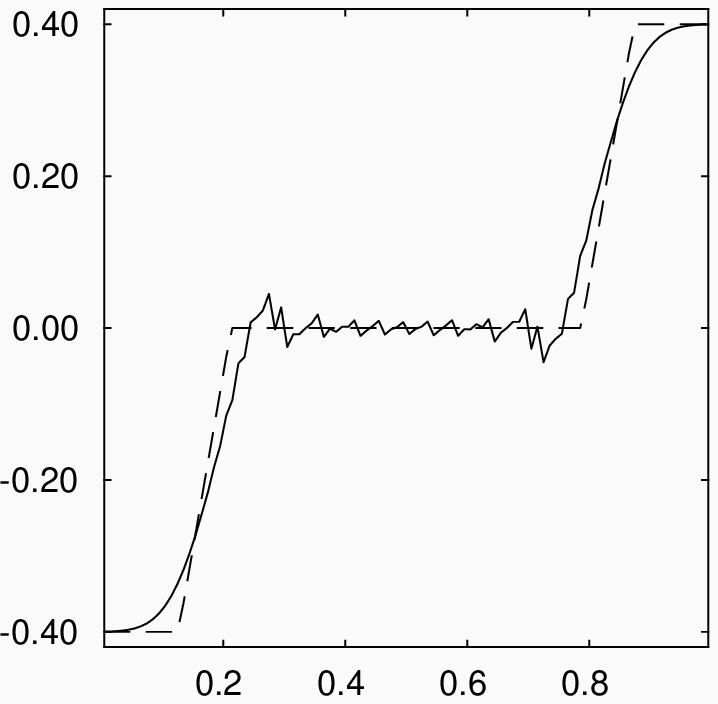}}
    \subfigure[Specific Entropy]{\includegraphics[scale = 0.63]{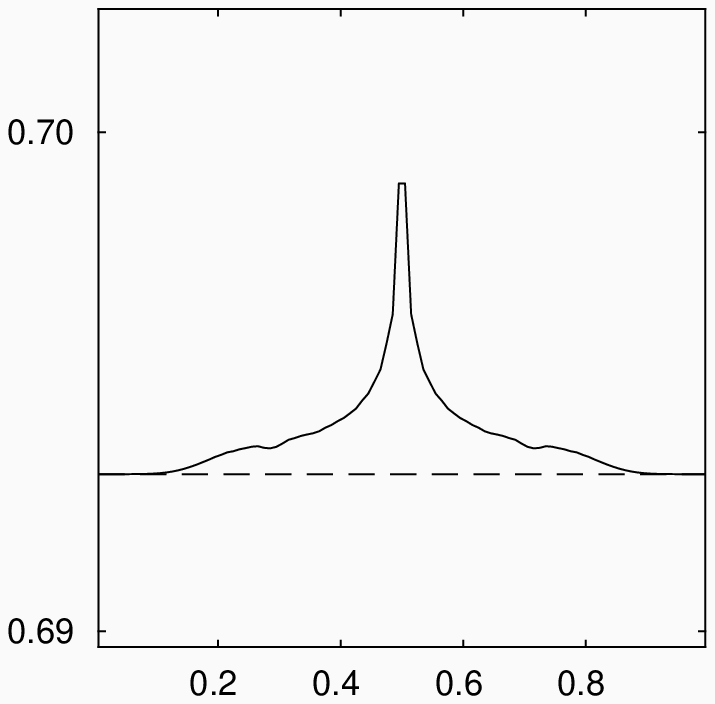}}
    \caption{Receding flow problem: Numerical solution (full line) at $t = 0.18s$ with the EC Roe flux and Backward Euler in time. 100 cells and $\Delta t = 10^{-3}s$. }
    \label{fig:RF_BEECRoe}
\end{figure}
\subsection{Fully discrete analysis with EC schemes}
The entropy rise observed with EC schemes can be explained if we consider the influence of time-integration. Let's assume that Forward Euler is used together with the EC Roe flux. We are interested in the jump of specific entropy $S^1 - S^0$ in the R cell after the first time step. We have the following discrete equation for density:
\begin{equation}\label{eq:density_fd}
    \rho^{1} - \rho^0 = \frac{\Delta t}{\Delta x} [(\rho u)_{1/2} - (\rho u)_{3/2}] = -\frac{\Delta t}{\Delta x} \rho^0 u^0. 
\end{equation}
Using the analysis of section 4, we have the following discrete equation for the ``entropy'' $U = - \frac{\rho S}{\gamma -1}$:
\begin{equation}\label{eq:entropy_fd_0}
    (\rho S)^{1} - (\rho S)^0 = \frac{\Delta t}{\Delta x} [(\rho u S)_{1/2} - (\rho u S)_{3/2}] + (1 - \gamma)\mathcal{E}^{FE}.
\end{equation}
The interface flux for the entropy is given by eq. (\ref{eq:entropy_interface}). With the EC Roe flux, we get the values:
\begin{equation*}
    (\rho u S)_{1/2} = 0, \ (\rho u S)_{3/2} = \rho^0 u^0 S^0.
\end{equation*}
The entropy production term $\mathcal{E}^{FE}$ is given by eq. (\ref{eq:TadFE}). Eq. (\ref{eq:entropy_fd_0}) becomes:
\begin{equation}\label{eq:entropy_fd}
    (\rho S)^{1} - (\rho S)^0 = -\frac{\Delta t}{\Delta x} \rho^0 u^0 S^0 + (1 - \gamma)\mathcal{E}^{FE}.
\end{equation}
Combining equations (\ref{eq:entropy_fd}) and (\ref{eq:density_fd}) gives:
\begin{equation*}
    (\rho S)^{1} - (\rho S)^0 = S^0 (\rho^1 - \rho^0) + (1 - \gamma)\mathcal{E}^{FE}.
\end{equation*}
Regrouping, one obtains:
\begin{equation}\label{eq:s_jump}
    S^1 - S^0 = (1 - \gamma)\mathcal{E}^{FE}/\rho^1.
\end{equation}
equation (\ref{eq:s_jump}) is exact and shows that when the EC Roe flux is combined with 1st order Euler explicit in time, the spurious entropy production is due to the entropy produced by the time-integration scheme. While this type of analysis can hardly be carried out for more complex combinations of fluxes and time-integration schemes, it is enough to justify the non-negligible impact of time-integration on the numerical behavior of the solution. Note that if non EC fluxes were used in space and if conservation of energy was replaced with conservation of entropy, the term $(1 - \gamma) \mathcal{E}^{FE}$ would disappear and we would indeed obtain that $S^1 - S^0 = 0$. The same result could be obtained with implicit schemes depending on the flux used in space.

\subsection{Numerical experiments}
\begin{figure}[h!]
    \centering
    \subfigure[Density]{\includegraphics[scale = 0.63]{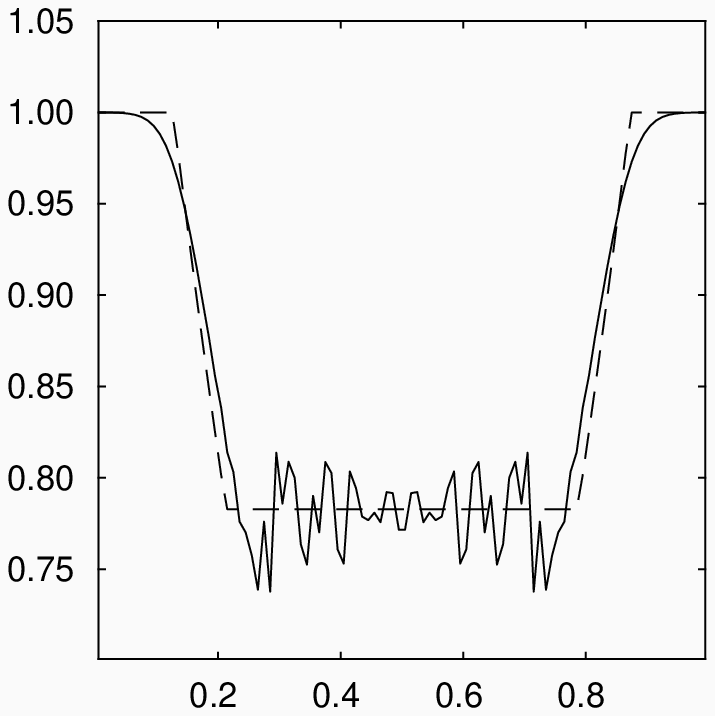}}
    \subfigure[Pressure]{\includegraphics[scale = 0.63]{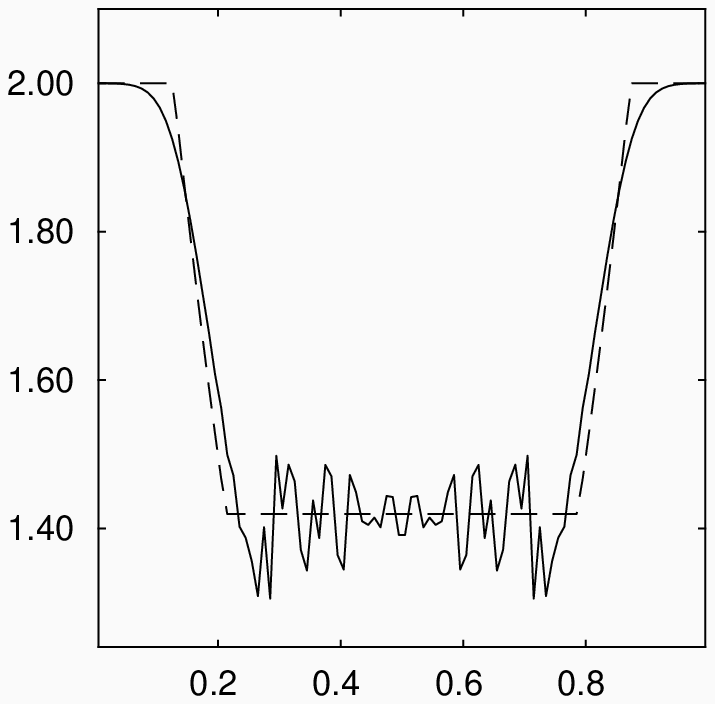}}
    \subfigure[Velocity]{\includegraphics[scale = 0.63]{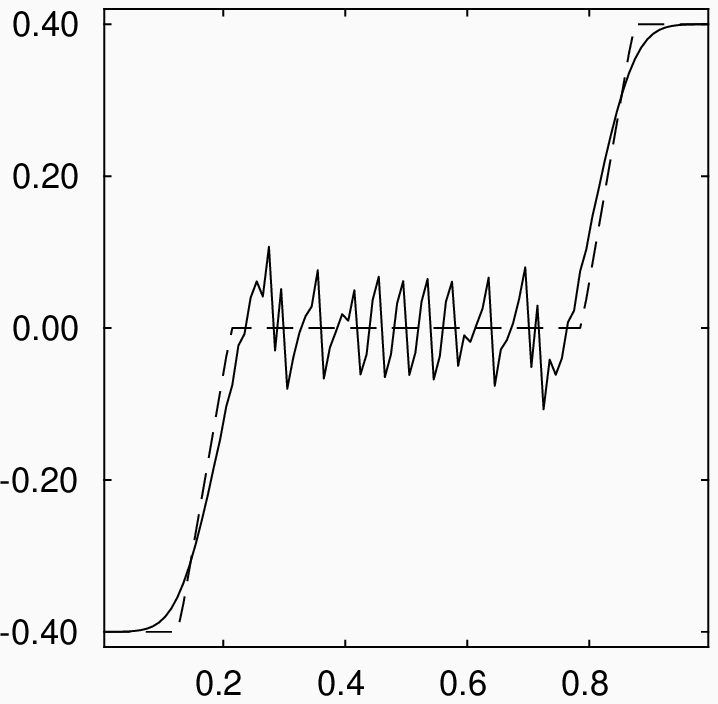}}
    \subfigure[Specific Entropy]{\includegraphics[scale = 0.63]{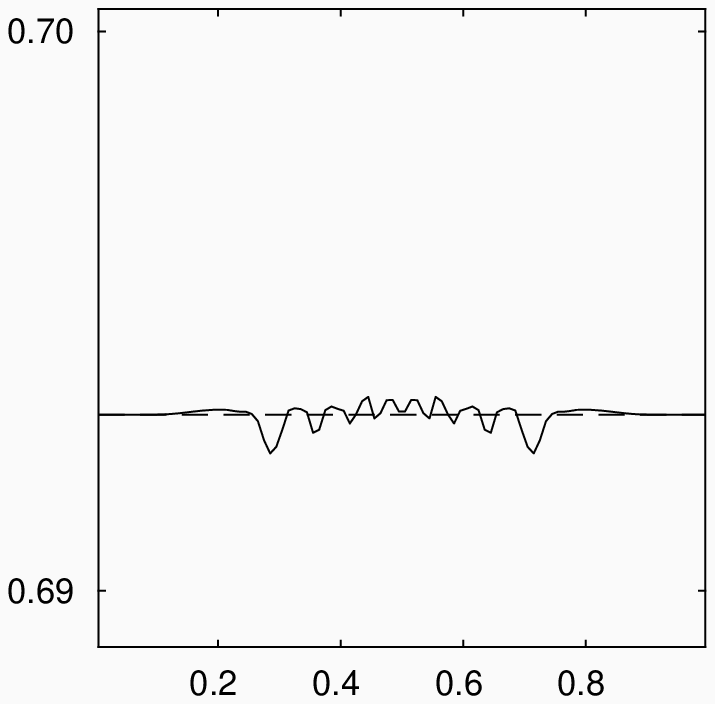}}
    \caption{Receding flow problem: Numerical solution (full line) at $t = 0.18s$ with the EC Roe flux and the EC scheme in time. 100 cells and $\Delta t = 10^{-3}s$. }
    \label{fig:RF_ECECRoe}
\end{figure}

\begin{figure}[h!]
    \centering
    \includegraphics[scale = 0.55]{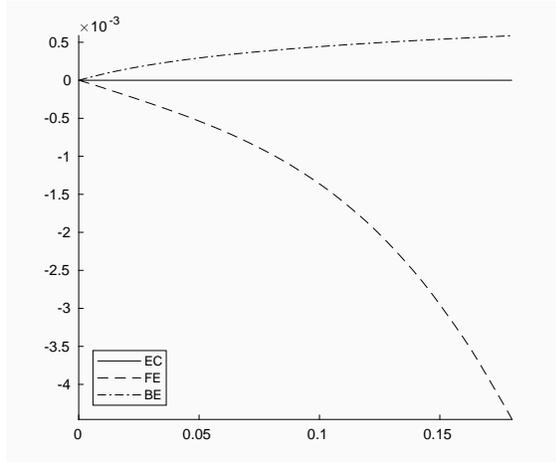}
    \caption{Receding flow problem: Total Entropy ($\rho S$) production over time. The EC Roe flux in used in space. EC: Our new EC scheme in time. FE: Forward Euler in time. BE: Backward Euler in time. 100 cells and $\Delta t = 10^{-3}s$. }
    \label{fig:RF_EntropyECRoe}
\end{figure}

\begin{figure}[h!]
    \centering
    \includegraphics[scale = 0.57]{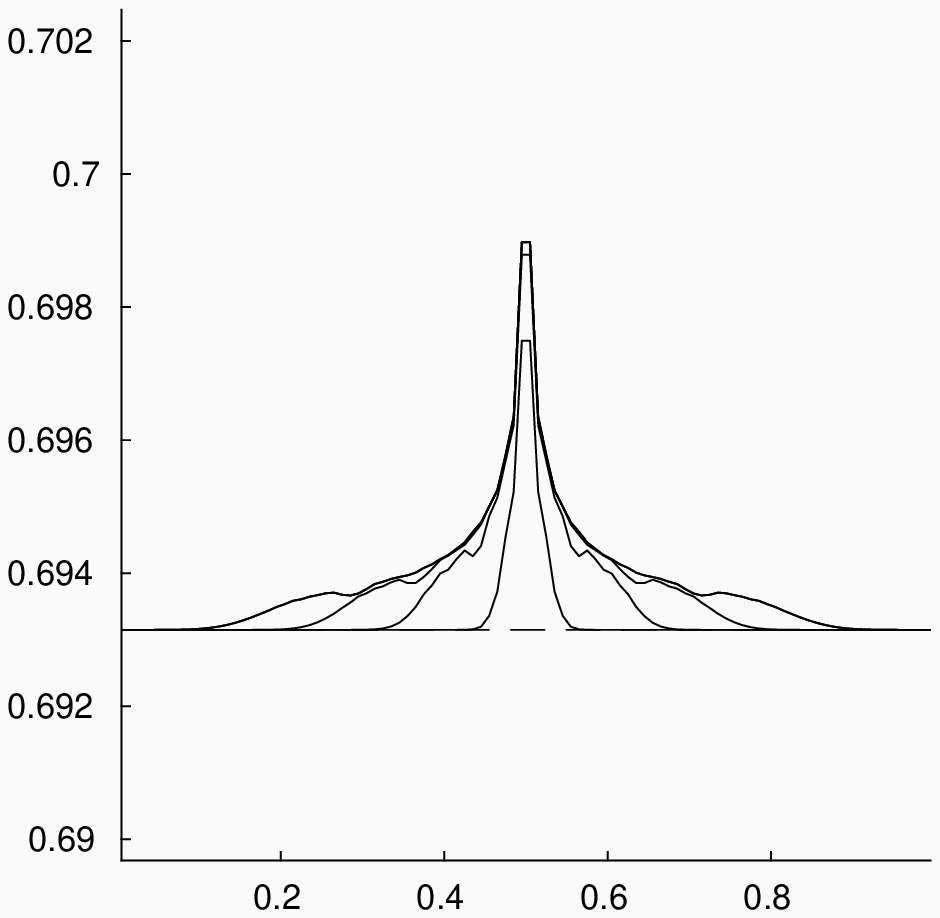}
    \caption{Receding flow problem: Snapshots of the specific entropy profile when the EC Roe flux is used in space and BE is used in time. 100 cells and $\Delta t = 10^{-3}s$. }
    \label{fig:RF_EC_Entropyevol_BE_100_1em3}
\end{figure}

In this section, we investigate at the numerical solution of the receding flow problem when an EC flux is used in space. Three time-integration schemes are considered: Forward Euler, Backward Euler and the EC time scheme we introduced in section 4.2. \\ 
\indent The EC time scheme and BE are implicit schemes. The nonlinear equations are solved iteratively using Newton's method followed by a line search algorithm. The Jacobian matrices associated with each scheme are directly computed using the complex step method. \\
\indent Figures \ref{fig:RF_FEECRoe}, \ref{fig:RF_BEECRoe} and \ref{fig:RF_ECECRoe} show the numerical solution with Forward Euler, Backward Euler, and the EC time scheme, respectively, for a grid of 100 elements and a time-step $\Delta t= 10^{-3}s$. Oscillations of various magnitudes are observed in all three cases. The oscillations are more pronounced in the Forward Euler case and their magnitude keeps growing with time because of the spurious entropy that is brought in by FE at every iteration. The opposite trend is observed in the Backward Euler case. The oscillations are much more damped. This is a consequence of the positive entropy production of BE. We also observe a specific entropy profile (figure \ref{fig:RF_BEECRoe}-(d)) that is similar to the overheated profile (figure \ref{fig:RF_classicRoe}-(d)) observed with conventional schemes. Figure \ref{fig:RF_EC_Entropyevol_BE_100_1em3} shows how the specific entropy profile in the BE case evolves over time. The spike at the center increases with time as the oscillations cause by the EC flux in space are damped by BE. In the case of a fully discrete conservative scheme (figure \ref{fig:RF_ECECRoe}), what we observe are oscillations of a magnitude higher than those of BE and much lower than for FE. The growth of oscillations is more controlled than for Forward Euler. We also note that the specific entropy profile is the most accurate of all three cases. The BE case still provides the best solution overall and illustrates why entropy-stable (as opposed to entropy conservative) schemes are used in practice. This shows that enforcing conservation of entropy, even when it is a property of the exact solution, does not necessarily lead to a better behaved numerical solution. Figure \ref{fig:RF_EntropyECRoe} shows the production of entropy over time in all three cases. It confirms the entropy stability properties of FE and BE that Tadmor proved, and the entropy conservation property of our time scheme. The rate at which BE produces positive entropy decreases with time. This is because, as time goes by, the oscillations caused by the EC flux in space are damped by the dissipation of BE and the numerical solution becomes smoother. 

\begin{figure}[h!]
    \centering
    \subfigure[Density]{\includegraphics[scale = 0.48]{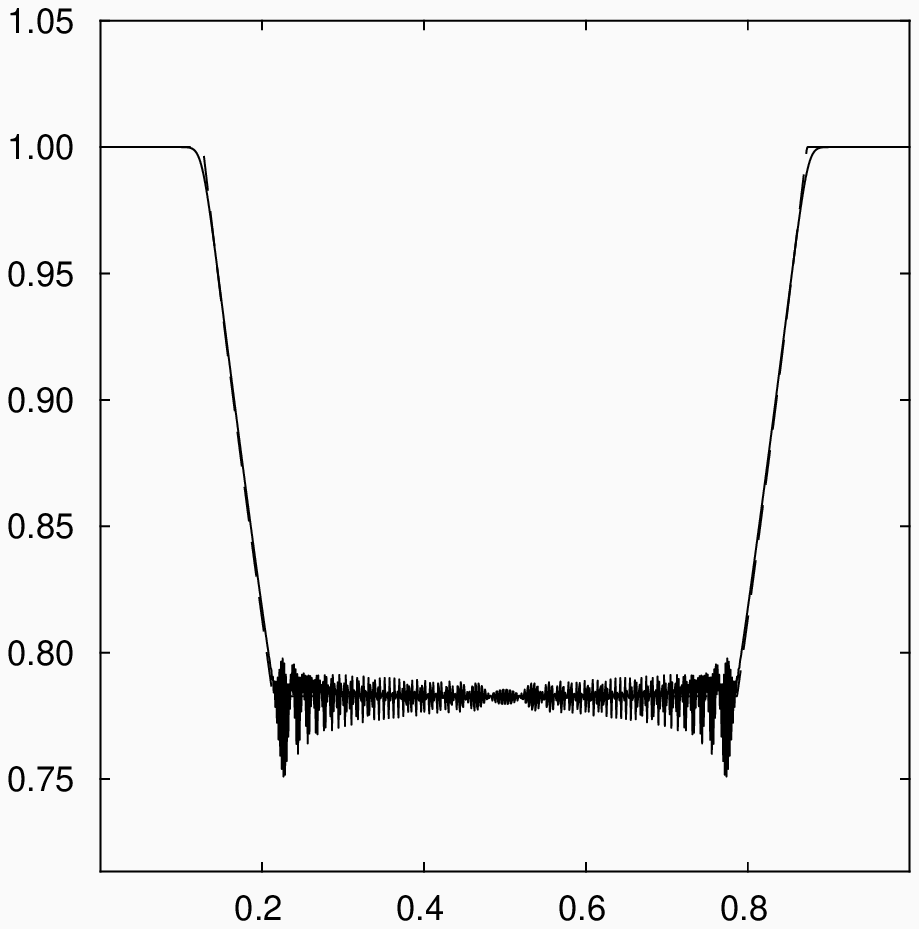}}
    \subfigure[Pressure]{\includegraphics[scale = 0.48]{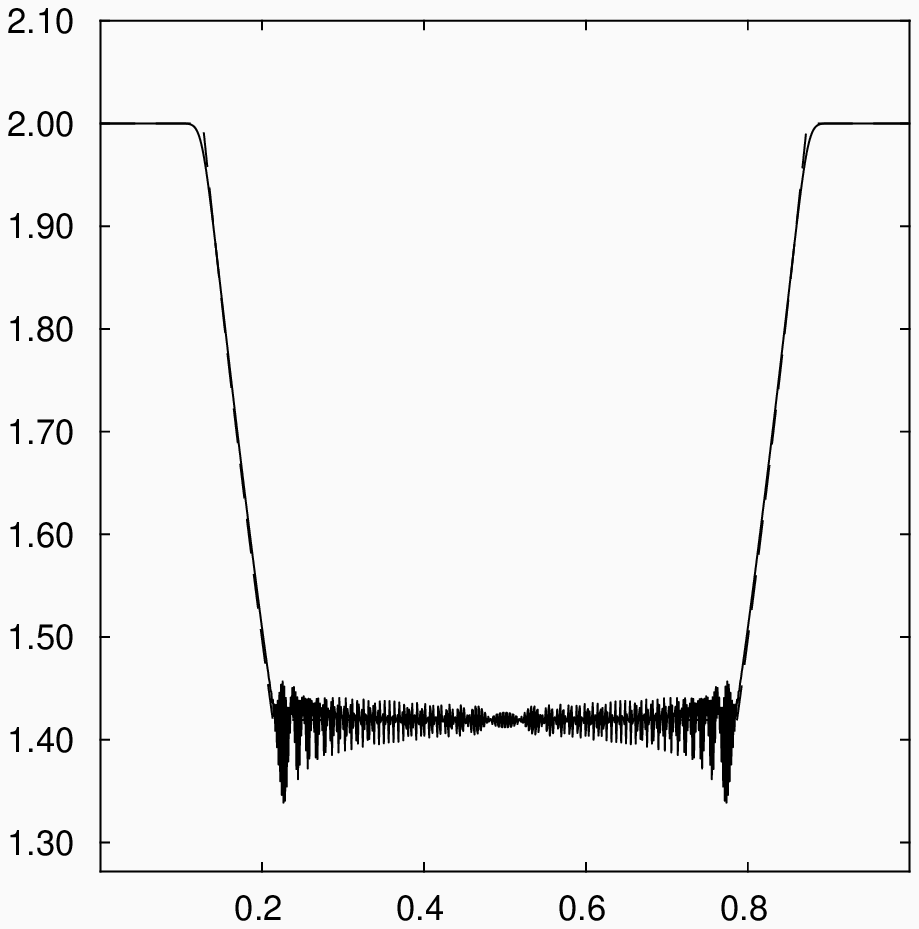}}
    \subfigure[Velocity]{\includegraphics[scale = 0.48]{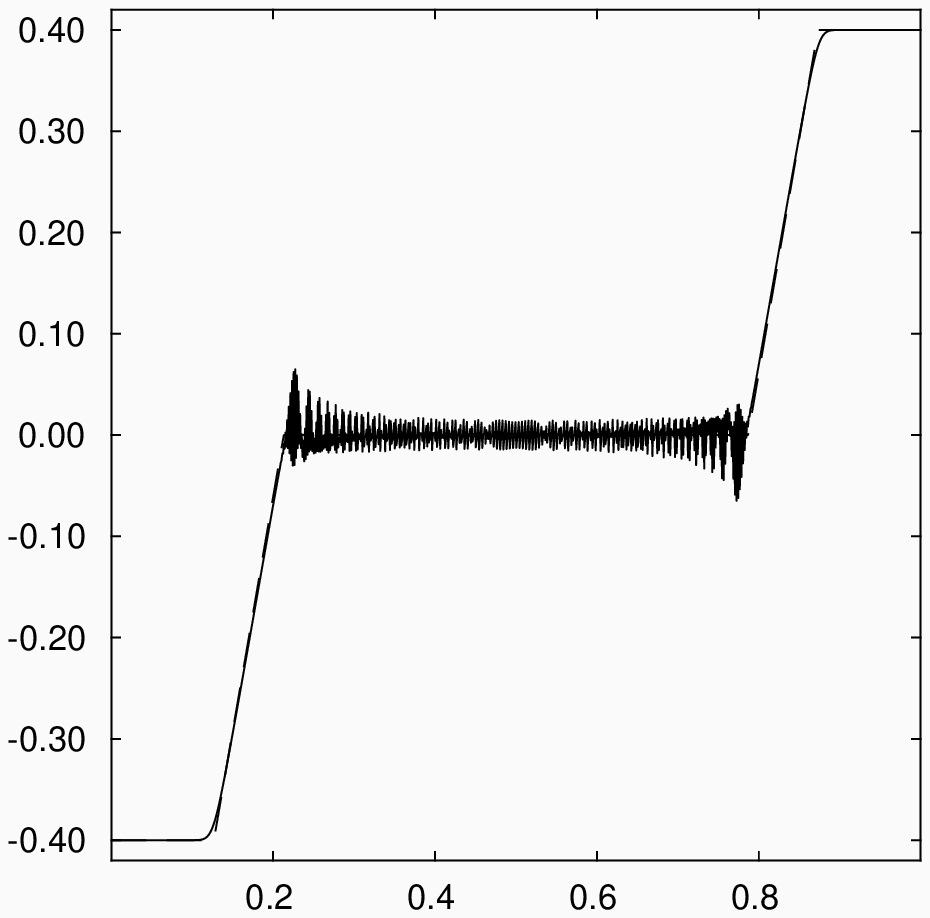}}
    \subfigure[Specific Entropy]{\includegraphics[scale = 0.48]{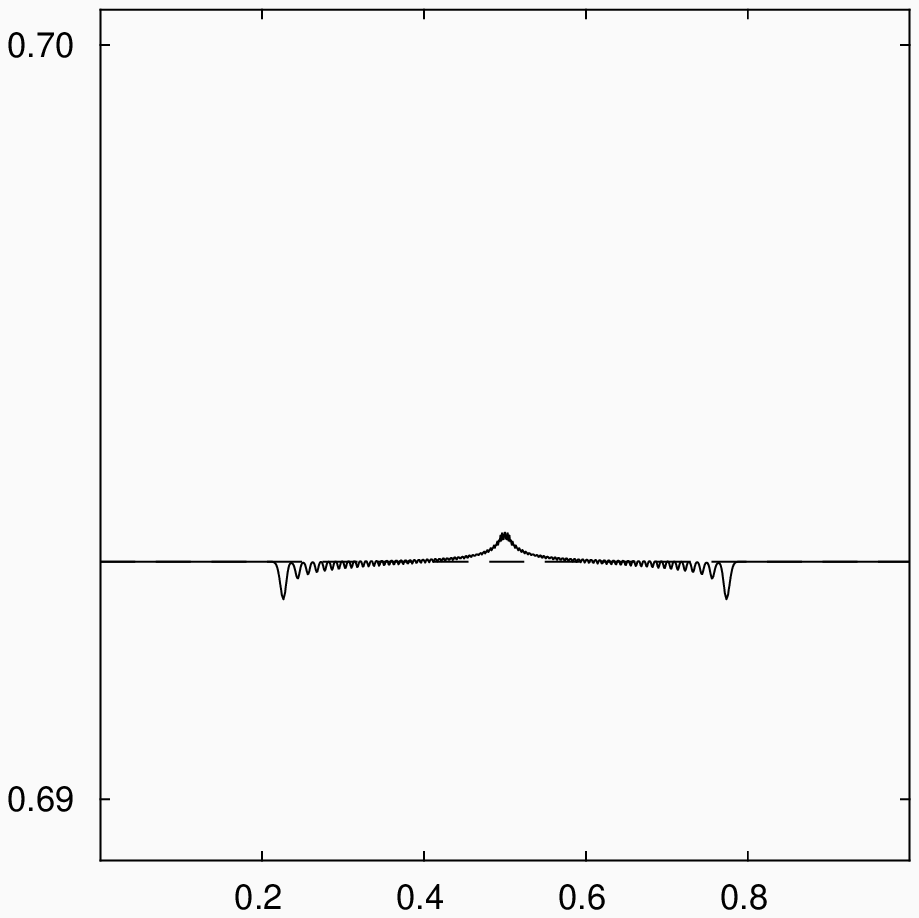}}
    \caption{Receding flow problem: Numerical solution (full line) at $t = 0.18s$ with the EC Roe flux and our EC scheme in time. 1000 cells and $\Delta t = 10^{-4}s$. }
    \label{fig:RF_ECECRoe_fine}
\end{figure}
\begin{figure}[h!]
    \centering
    \includegraphics[scale = 0.57]{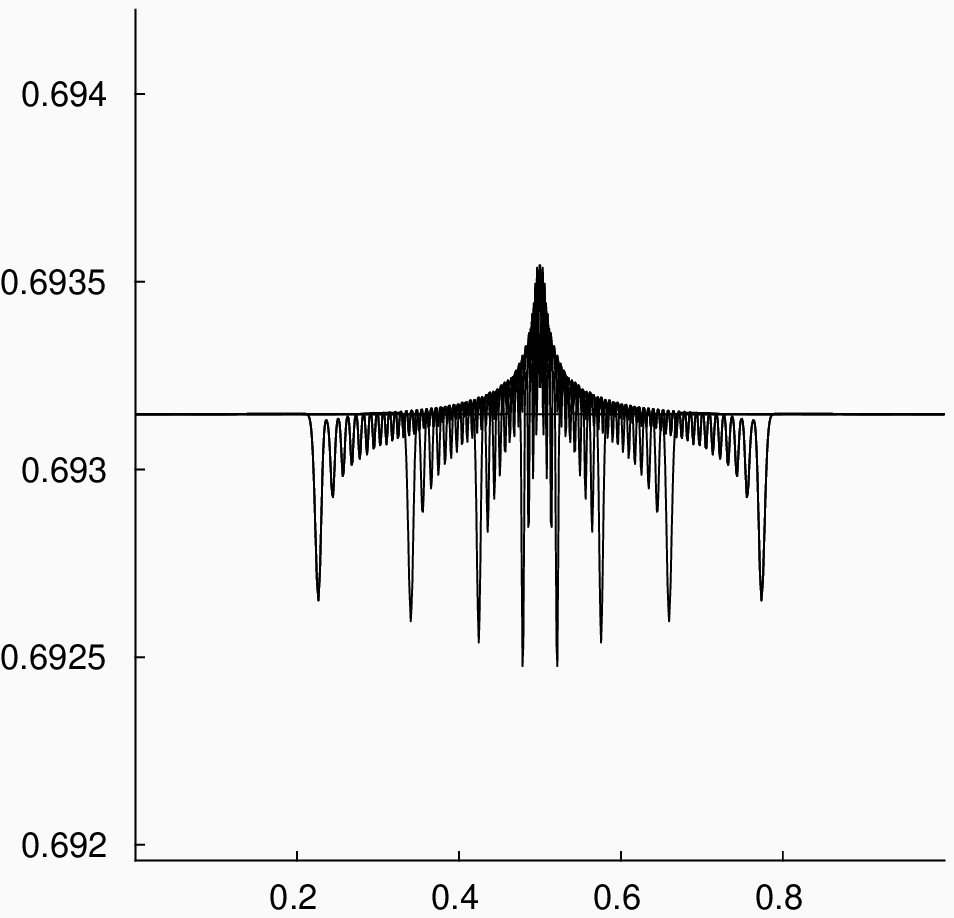}
    \caption{Receding flow problem: Snapshots of the specific entropy profile when the EC Roe flux is used in space and our EC scheme is used in time. 1000 cells and $\Delta t = 10^{-4}s$. }
    \label{fig:RF_EC_Entropyevol}
\end{figure}
\begin{figure}[h!]
    \centering
    \subfigure[Density]{\includegraphics[scale = 0.48]{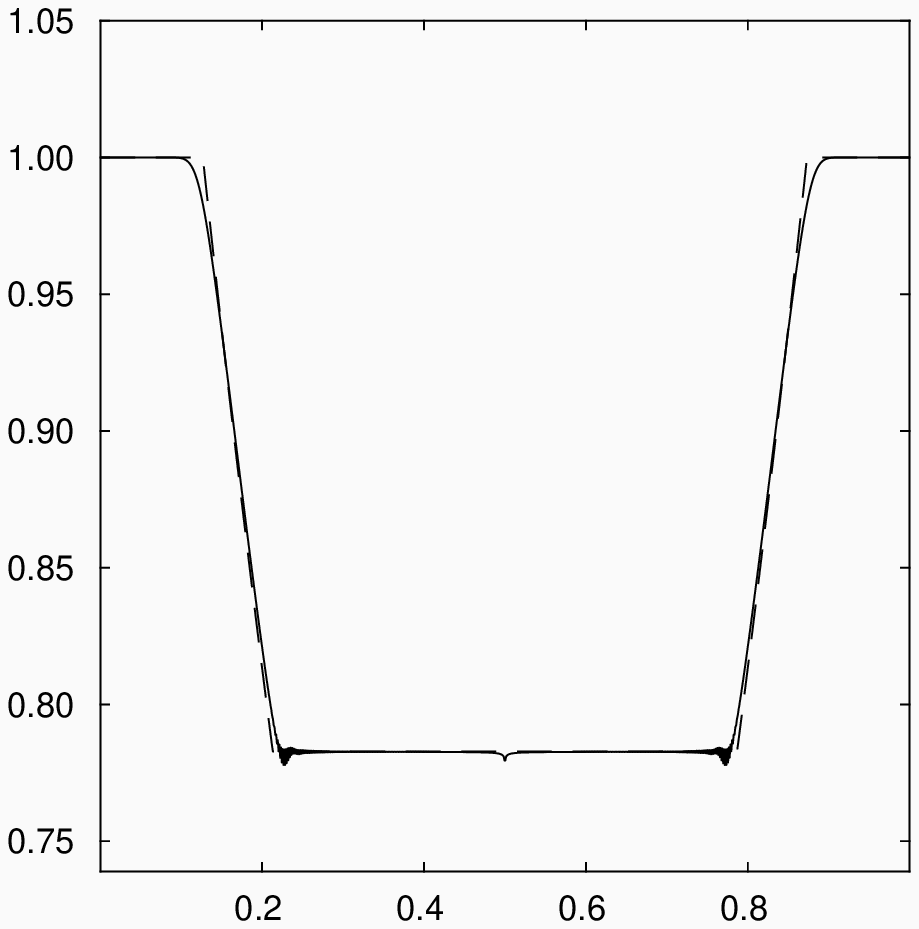}}
    \subfigure[Pressure]{\includegraphics[scale = 0.48]{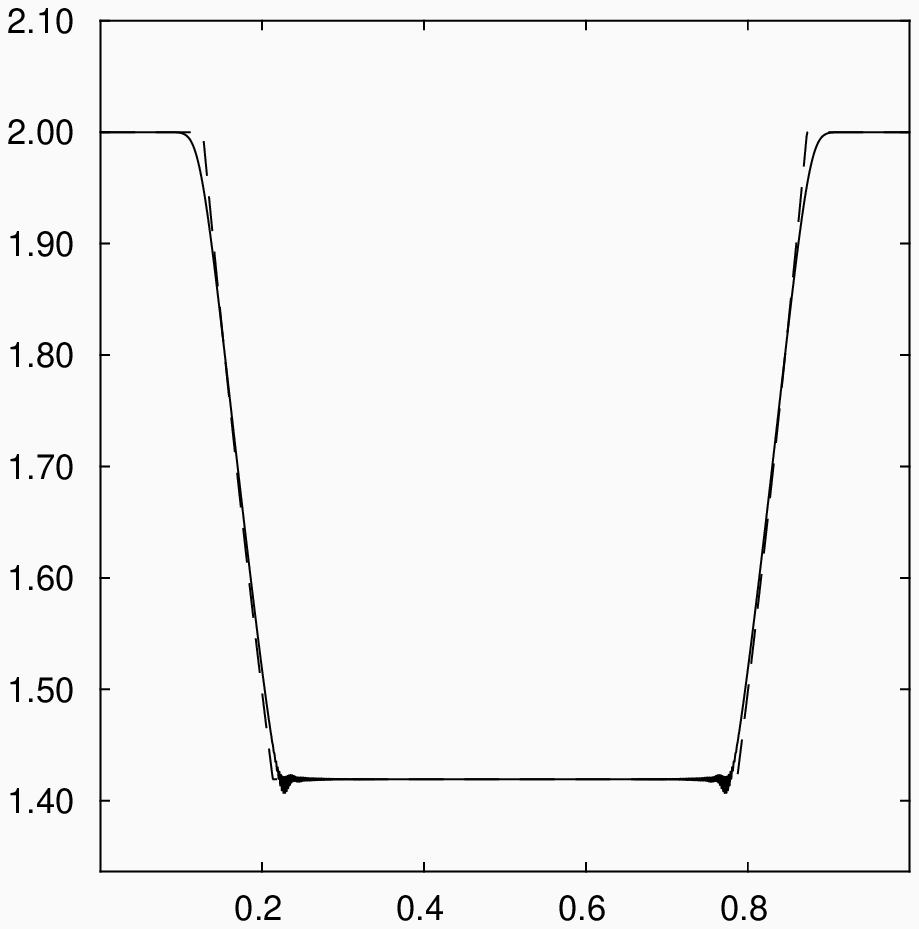}}
    \subfigure[Velocity]{\includegraphics[scale = 0.48]{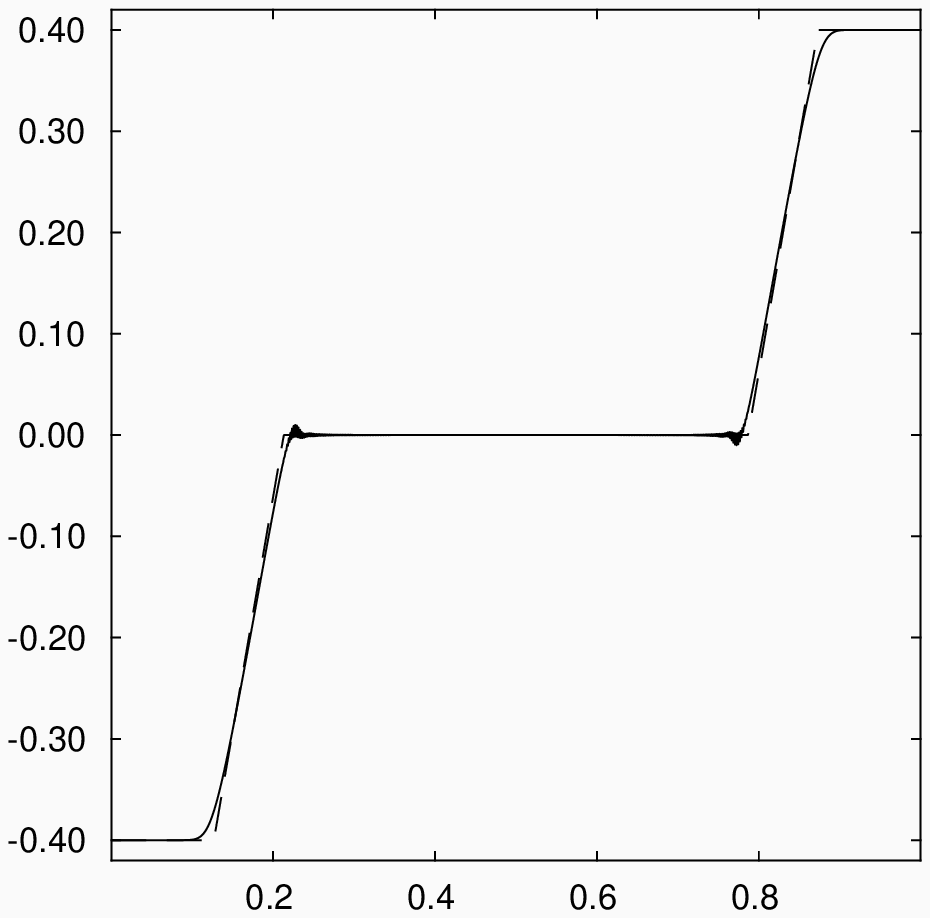}}
    \subfigure[Specific Entropy]{\includegraphics[scale = 0.48]{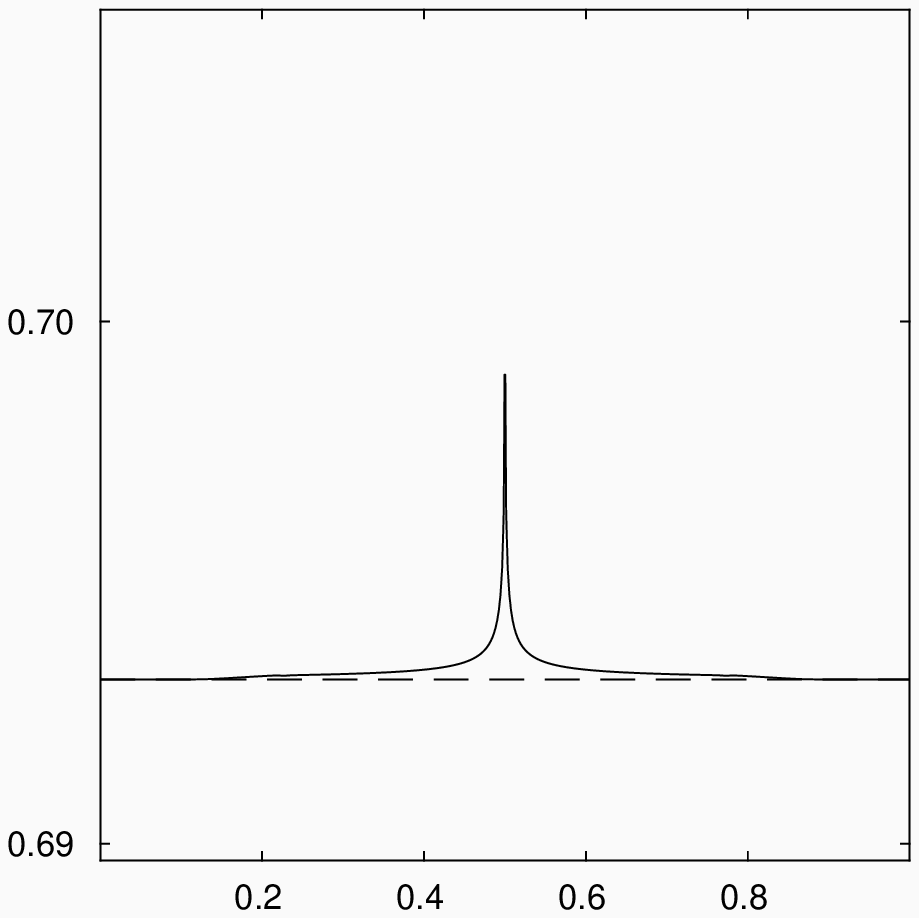}}
    \caption{Receding flow problem: Numerical solution (full line) at $t = 0.18s$ with the EC Roe flux and Backward Euler in time. 1000 cells and $\Delta t = 10^{-4}s$. }
    \label{fig:RF_BEECRoe_fine}
\end{figure}

Figure \ref{fig:RF_BEECRoe} suggests that when an EC flux is used in space, the overheating is correlated to a positive entropy production. Increasing the time step for BE increases that production. Taking $\Delta t = 4.3 \times 10^{-3}s$ results in a solution that is very close to the obtained with FE in time and the classic Roe flux in space (figure \ref{fig:RF_classicRoe}). More generally, it appears from figures \ref{fig:RF_FEECRoe}-(d), \ref{fig:RF_BEECRoe}-(d) and \ref{fig:RF_ECECRoe}-(d) that the entropy production of the fully discrete scheme has a strong impact on the specific entropy profile that Liou's latest studies focused on. \\
\indent Figures \ref{fig:RF_ECECRoe_fine} and \ref{fig:RF_BEECRoe_fine} show the numerical solution with a finer mesh (1000 cells) and a smaller time-step ($\Delta t = 10^{-4}s$) for the EC time scheme and BE. FE does not converge for this configuration. In figure \ref{fig:RF_ECECRoe_fine}, the oscillations have a higher frequency and a lower magnitude than in the coarser configuration. The oscillations obtained with the EC scheme in time appear to have more fine-grained structures than in figure \ref{fig:RF_ECECRoe}. The envelope of the density oscillations observed in figure \ref{fig:RF_ECECRoe_fine}-(a) and the density profile of the BE case \ref{fig:RF_BEECRoe_fine}-(a) suggests that the fully discrete EC scheme essentially under-estimates density before the rarefaction waves. We can see that the specific entropy profile (figure \ref{fig:RF_ECECRoe_fine}) produced by the fully discrete EC scheme has the spike associated with the overheating at its center, and two spikes at the beginning of each rarefaction that seem to compensate for the overheating. Figure \ref{fig:RF_EC_Entropyevol}, which features snapshots of the specific entropy profile over time, shows that this structure is conserved over time. In addition, the specific entropy increases with time everywhere in the domain.   
\indent 
\subsection{Additional remarks}

\indent The EC scheme that has been developed here is just one way among many to conserve entropy in addition to mass, momentum and energy at the fully discrete level. In the scalar case $(N = 1)$, the entropy conservation condition in space (\ref{eq:ECcond2}) has only one solution. For systems $(N > 1)$ there is more than one possible EC flux. Likewise, the intermediate state in time we used in our EC time scheme is just one choice among many. The time scheme given by eq. (\ref{eq:LeFloch}) is part of a more general class of entropy conservative scheme that LeFloch \textit{et al.} introduced in \cite{LeFloch}. While one can arguably take the stance that fully discrete EC schemes will produce a similar behavior to that in figures \ref{fig:RF_ECECRoe} and \ref{fig:RF_ECECRoe_fine}, it is known from past work \cite{Chandra, Roe, Derigs} that all EC fluxes do not perform equally. Besides simplicity, one of the reasons why the EC Roe flux is preferred over the first EC flux of Tadmor is that the latter does not preserve stationary contact discontinuities. Chandrasekhar \cite{Chandra} introduced an EC flux that has the additional property of discretely preserving, in the sense of Jameson \cite{Jameson}, the kinetic energy of the system. This type of property is often sought when turbulent flows are simulated. Another metric is how good of a foundation an EC flux constitutes in an ES scheme. The dissipation component of entropy stable fluxes is often seen as the complement needed by EC fluxes in the presence of shocks. An EC scheme will produce nonphysical solutions (oscillations) in the presence of shocks because entropy is not produced. This picture is correct but incomplete. In the presence of rarefaction waves and moving contact discontinuities, which do not physically require any production of entropy, EC schemes have the same oscillatory behavior (the receding flow problem is an illustration). This places an additional burden on the dissipation term which has to make up for the flaws of its foundation. Derigs \textit{et al.} \cite{Derigs} showed that  entropy stable schemes perform better on high-pressure shock problems if Chandrasekhar's EC flux is used instead of Roe's. \\
\indent Besides, one must not forget Liou's finding that the overheating occurs in the very first instants where the two receding rarefactions waves introduce the most significant discontinuities. It can be easily checked that with a fine enough grid, any scheme, EC or non-EC, will perform well if started from the exact solution at $t = \tau > 0$. Completely excluding EC schemes would require showing that no fully discrete EC scheme can survive in that $[0, \ \tau]$ time window. We are not able to prove or disprove such a claim. 
\section{Conclusions}
\label{sec:conclusions}
In this work, entropy conservative schemes, which allow for the conservation of an additional quantity at the semi-discrete level, were considered for the receding flow problem. This was motivated by Liou's latest study that showed the connection between the overheating and a spurious entropy production \textit{ab initio}. While Liou's semi-discrete analysis suggested that the EC flux of Roe would prevent the overheating, a fully discrete analysis showed the influence of time-integration on the entropy production. \\
\indent This observation brought about the question of the behavior of a fully discrete entropy conservative scheme on this type of problem. Building on the analogy between the entropy conservation condition for the spatial fluxes and the entropy conservation condition of a class of time-integration schemes considered by LeFloch \textit{et. al}, we derived a new entropy conservative time-integration scheme. Combining it with an EC flux, we observed that it does not necessarily make the solution better. A better specific entropy profile is obtained but the oscillatory nature of the numerical solution does not make it a practical option.\\
\indent Whether all entropy conservative discretizations would have the same unsatisfactory behavior on this type of problem, where one expects the continuous solution to conserve entropy, is a question that requires further investigation.  

\section*{Acknowledgments}
Ayoub Gouasmi would like to thank Laslo Diosady and Nicholas Burgess for insightful conversations from which this work benefited. The authors also thank Professor Philip Roe for helpful comments on the manuscript.

%\Appendix
%\section{The use of appendices}
%The \verb|\appendix| command may be used before the final sections
%of a paper to designate them as appendices. Once \verb|\appendix|
%is called, all subsequent sections will appear as 

%\appendix
%\section{Title of appendix} Each one will be sequentially lettered
%instead of numbered. Theorem-like environments, subsections,
%and equations will also have the section number changed to a letter.

%If there is only {\em one} appendix, however, the \verb|\Appendix|
%(with a capital letter) should be used instead. This produces only 
%the word {\bf Appendix} in the section title, and does not add a letter.
%equation numbers, theorem numbers and subsections of the appendix
%will have the letter ``A''  designating the section number.

%If you don't want to title your appendix, and just call it
%{\bf Appendix A.} for example, use \verb|\appendix\section*{}| 
%and don't include anything in the title field. This works
%opposite to the way \verb|\section*| usually works, by including the
%section number, but not using a title.

%Appendices should appear before the bibliography section, not after,
%and any acknowledgments should be placed after the appendices and before
%the bibliography. 

\bibliographystyle{siamplain}

\end{document}